\documentclass[12pt,a4paper]{article}

\usepackage[dvipdfmx]{color}

\usepackage[pdftex]{graphicx}

\usepackage{authblk}

\usepackage{amsmath}
\usepackage{amsthm}
\usepackage{ascmac}
\usepackage{amssymb}
\usepackage{amsfonts}
\usepackage{mathtools}
\usepackage{bm}
\usepackage{here}
\usepackage{cancel}
\usepackage{comment}
\usepackage{fancyhdr}
\usepackage{cases}
\usepackage{multicol}
\usepackage{braket}
\usepackage{geometry}
\theoremstyle{definition}

\newtheorem*{theorem*}{Theorem}

\newtheorem*{definition*}{Definition}

\newtheorem*{remark*}{Remark}

\newtheorem*{proposition*}{Proposition}

\title{Representation of Geometric Objects by Color}
\author{Yusuke Imai \thanks{CONTACT: 93imaiyusuke@gmail.com}}
\affil{Graduate School of Engineering Science, Osaka University, Toyonaka, Osaka 560-8531, Japan}
\date{\today}

\begin{document}
\maketitle

\begin{abstract}
By introducing various actions involving color to geometrical objects, we represent a cube and simplex in four or fewer dimensions, the geometrical net of a cube and simplex in five or fewer dimensions, hyperprisms, truncated polytopes, stellated polytopes, and fractals such as Cantor dust and Menger sponge, and propose the ``four-dimensional" Menger sponge whose volume is infinite. 
\end{abstract}
 
 \tableofcontents

\newpage
\section{Introduction}\label{Introduction}
Color is often used to represent 3-dimensional data by 2-dimensional geometric objects. A typical example is the color representation of temperature in weather reports. In addition, complex 3-dimensional data in academic research is often represented by a 2-dimensional color map because of visibility and efficiency.

In this paper, we propose ways to represent geometrical objects by color. Such an attempt would have been partially made in a variety of settings (For example, see Ref.\cite{Irons}). This paper presents the unified methodology by using various actions on geometrical objects: coloring, uncoloring, $1/n$-coloring.

\begin{figure}[H]
\centering
{%
\resizebox*{2cm}{!}{\includegraphics{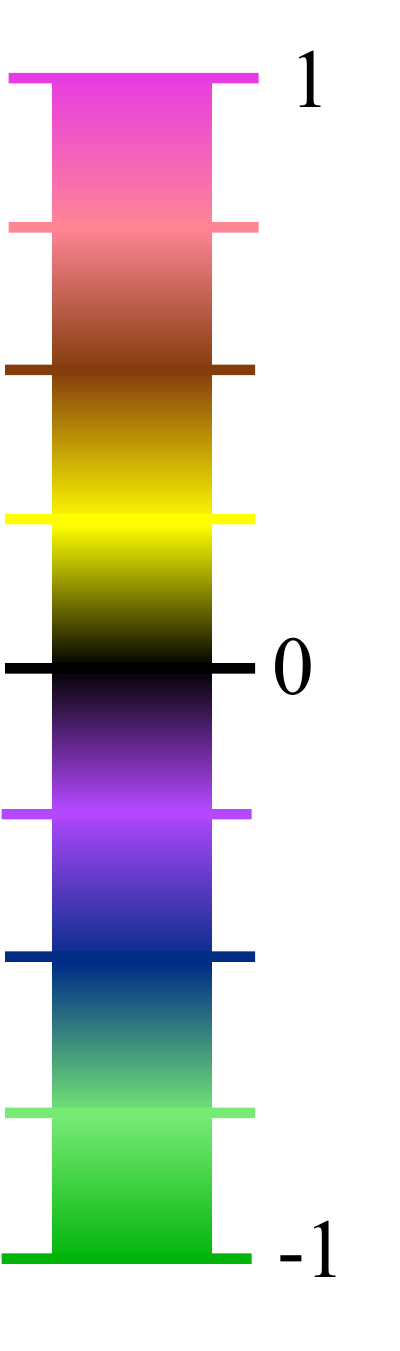}}}\hspace{5pt}
\caption{The color bar.} \label{color_bar}
\end{figure}

\begin{figure}[H]
\centering
{%
\resizebox*{12cm}{!}{\includegraphics{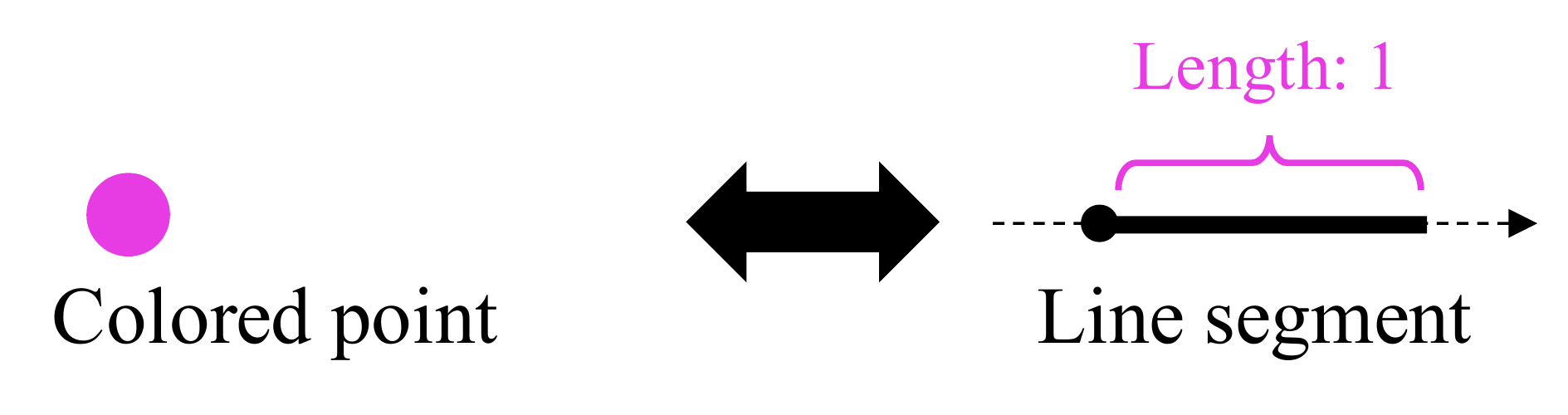}}}\hspace{5pt}
\caption{The colored point.} \label{color_point}
\end{figure}

\begin{figure}[H]
\centering
{%
\resizebox*{\textwidth}{!}{\includegraphics{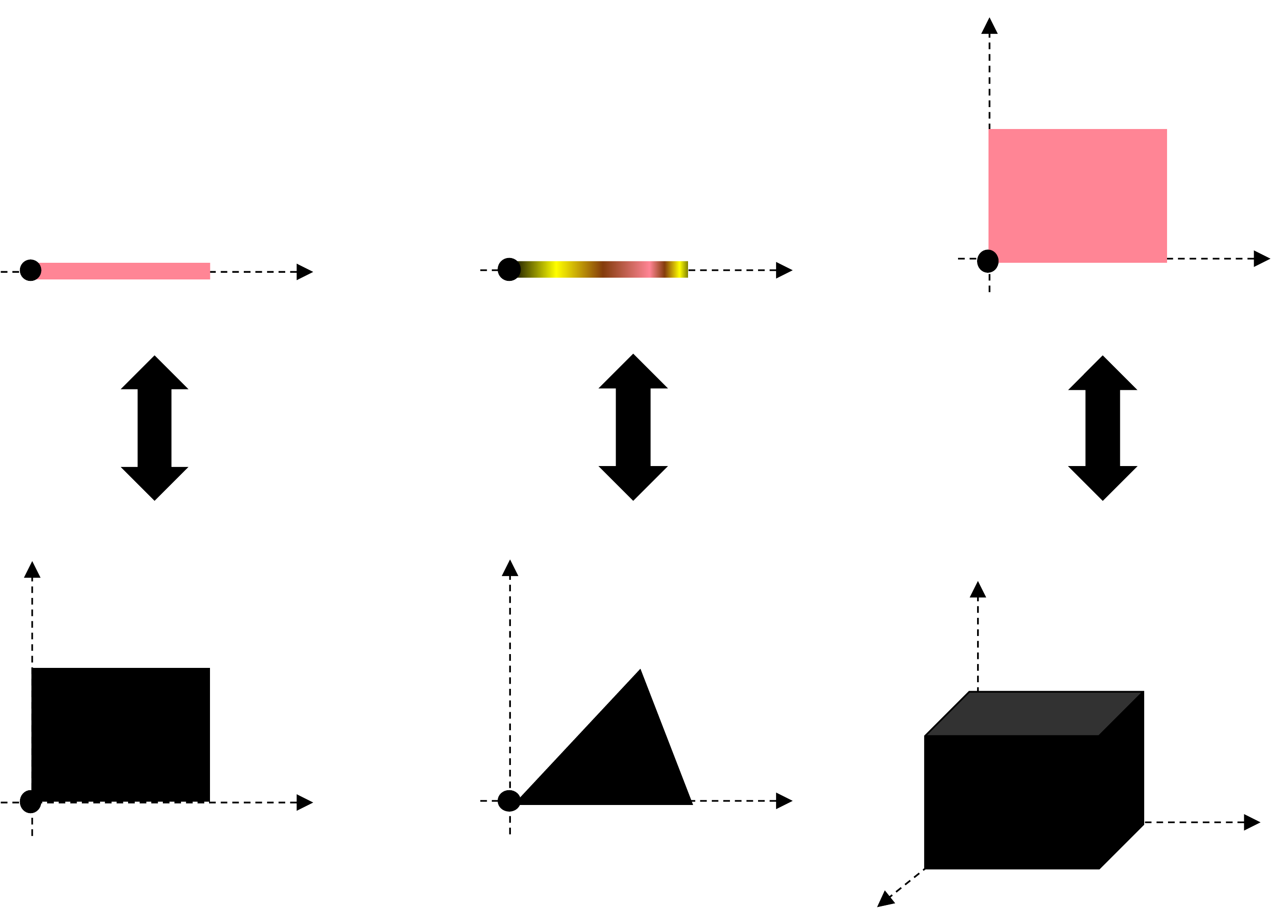}}}\hspace{5pt}
\caption{Examples of the color representation of a geometrical object: rectangle, triangle, rectangular.} \label{example}
\end{figure}

In this Introduction, first, we introduce color point. The colored point represents a line segment that exists in $\mathbb{R}^{1}_+ \coloneqq \{x_1 \, | \, x_1 \geq 0\}$, one of whose boundary is set on the color point, and whose length is determined by kinds of color as shown in Fig. \ref{color_bar}. The colored point is schematically shown in Fig. \ref{color_point}. The colored point is naturally generalized to higher dimensions such as colored line segment, colored plane, etc. For example, a rectangle/triangle/rectangular is represented by using the colored line segment as shown in Fig. \ref{example}.

In Sec. \ref{Regular Polytope}, we deal with $n$-cube and $n$-simplex. Especially, by using the color point and its generalizations, we propose a way to represent 4-cube and 4-simplex and the geometrical net of 5-cube and 5-simplex \cite{Schlafli,Coxeter} systematically. The color representation of the geometrical net of 4-cube, 4-simplex, 5-cube, and the 5-simplex leads to an intuitive way to answer questions such as ``How many edges are there in 5-cube?" or ``How many faces are there in 5-simplex?". As for the color representation of the geometrical net of $n$-simplex, we also use $A$-folding which is an operation of constructing the polytope $A$ from the geometrical net of $A$.

In Sec. \ref{Hyperprism}, we apply the coloring patterns that appeared in the color representation of $n$-cube and $n$-simplex to another class of polytopes, prism. First, we propose a way to represent polytopal prisms, i.e., polytopes represented by the Cartesian product of a polytope and a line segment. Then, we also propose a way to represent duoprisms, i.e., the Cartesian product of two polygons.

\begin{figure}[H]
\centering
{%
\resizebox*{12cm}{!}{\includegraphics{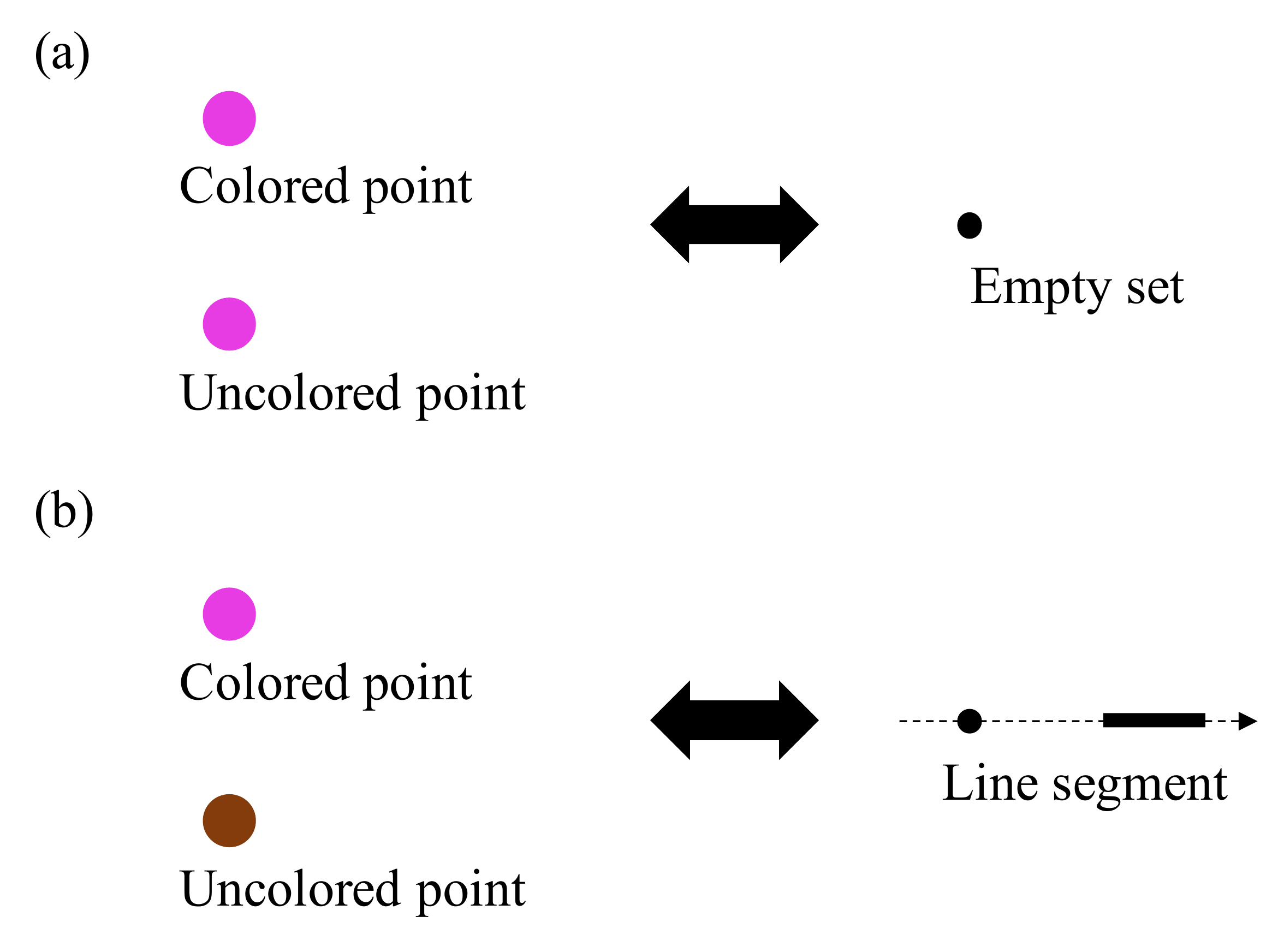}}}\hspace{5pt}
\caption{Examples of the coloring and uncoloring. (a) The pink-colored point and the pink-uncolored point make an empty set. (b) The pink-colored point and the brown-uncolored point make a line segment of length 0.5.} \label{uncolor}
\end{figure}

In Sec. \ref{Uncoloring and Truncation}, we use the concept of uncoloring. The uncolored point erases a part of a line segment represented by a colored point and its rate is also given by color. For example, a pink point represents a line segment whose length is 1, and a pink uncolored point erases it completely (see Fig. \ref{uncolor}(a)). Also, a pink point represents a line segment whose length is 1, and a brown uncolored point erases it in half (see Fig. \ref{uncolor}(b)). Then, we generalize the uncolored point to higher dimensions naturally. By using the uncolored point and its generalizations, we propose a way to represent truncated polytopes, i.e., polytopes obtained by cutting polytope vertices.

\begin{figure}[H]
\centering
{%
\resizebox*{\textwidth}{!}{\includegraphics{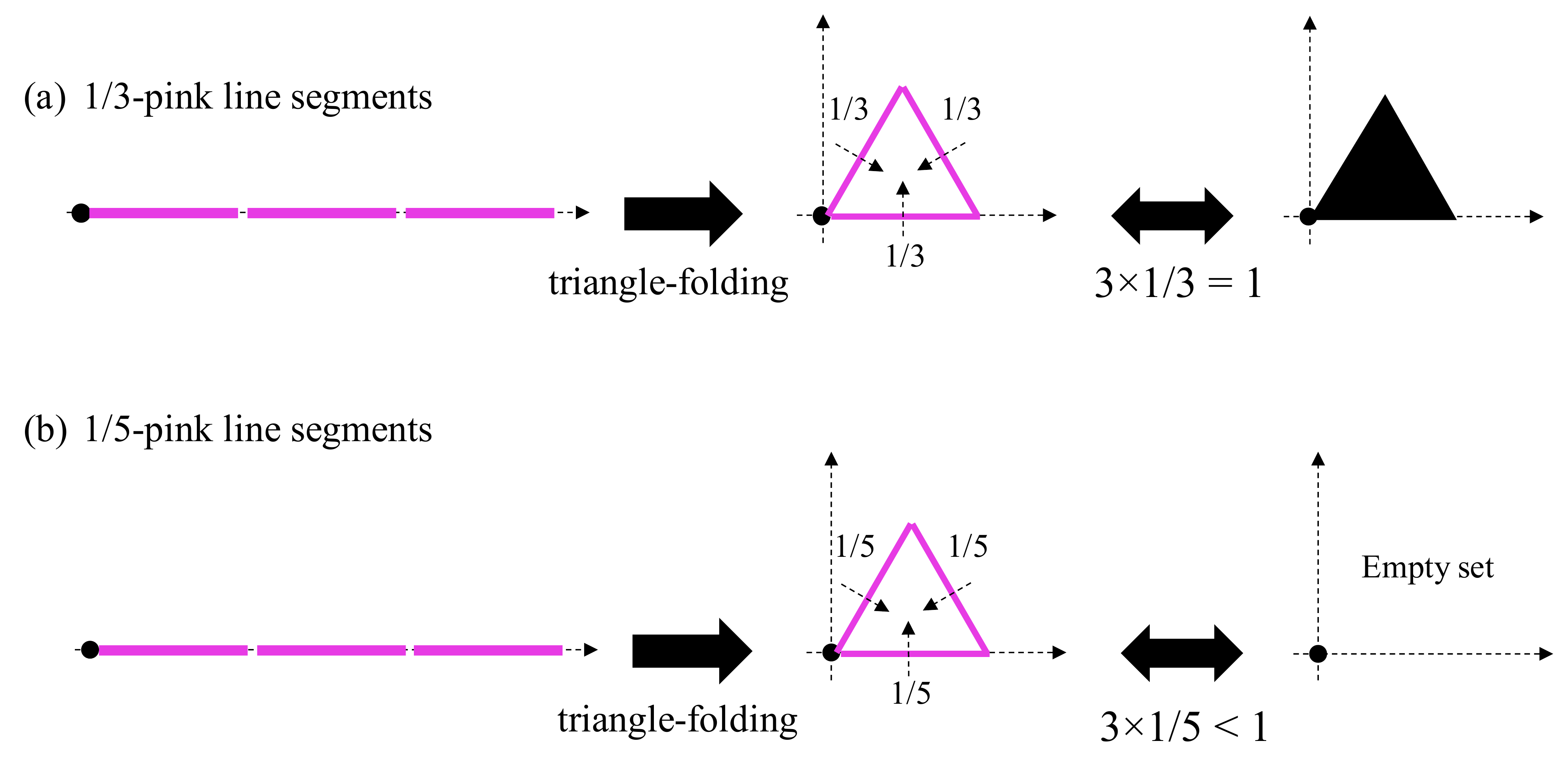}}}\hspace{5pt}
\caption{Examples of the $1/n$-coloring. Triangle-folding of successive (a) $1/3$-pink line segments and (b) $1/5$-pink line segments.}  \label{1/n}
\end{figure}

In Sec. \ref{1/n-Coloring and Stellation}, we use $1/n$-coloring. The $1/n$-colored figure does not represent any figures by itself. The $1/n$-coloring of a figure $A$ represents $1/n$-figure that has the same shape as the figure represented by the color representation of $A$ and we consider $n$-fold overlap of the $1/n$-figures to be identical to usual the color representation. For example, the triangle-folding of three successive $1/3$-pink line segments represents an equilateral triangle (see Fig.\ref{1/n}(a)) while the triangle-folding of three successive $1/5$-pink line segments does not represent any figures (see Fig.\ref{1/n}(b)). Then, by using $1/n$-colored line segments, we propose a way to represent star polygons.

In Sec. \ref{Fractal}, we represent two-dimensional Cantor dust (three-dimensional Cantor dust)\cite{Irons} by the coloring and uncolroing of Cantor set (two-dimensional Cantor dust), and Menger sponge\cite{Irons} by the coloring and uncolroing of Sierpinski carpet\cite{Irons}, and propose the ``four-dimensional" Menger sponge where each iteration is defined in the four-dimensional Euclidean space and whose volume is infinite.




\newpage
\section{Regular Polytope}\label{Regular Polytope}
\subsection{Hypercube}\label{Hypercube}
A one-dimensional hypercube (1-cube) is a line segment. A unit line segment is represented by a pink point (Fig.~\ref{color_point}).

\begin{figure}[H]
\centering
{%
\resizebox*{\textwidth}{!}{\includegraphics{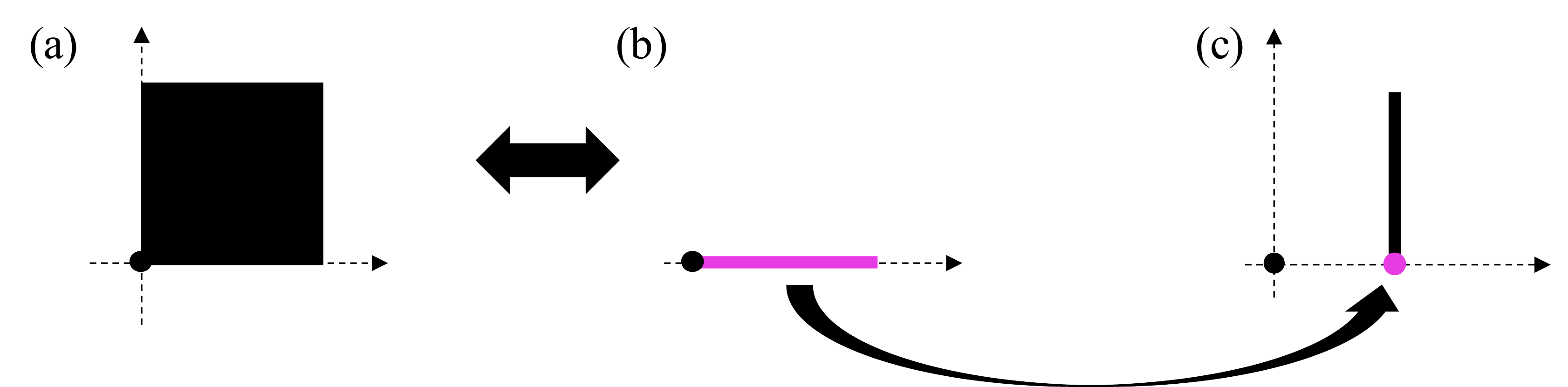}}}\hspace{5pt}
\caption{(a) 2-cube. (b) Color representation of 2-cube. (c) Decomposition of the color representation of 2-cube.}  \label{cube_2d}
\end{figure}

A two-dimensional hypercube (2-cube) is a square. A unit square (Fig.~\ref{cube_2d}(a)) is represented by a pink line segment whose length is 1 (Fig.~\ref{cube_2d}(b))\footnote{We call this uniform coloring square coloring.}. This colored figure can be interpreted as an infinite set of pink points corresponding to the unit line segments obtained by extending the pink points in the direction of the 2nd dimension as shown in Fig.~\ref{cube_2d}(c).

\begin{figure}[H]
\centering
{%
\resizebox*{10cm}{!}{\includegraphics{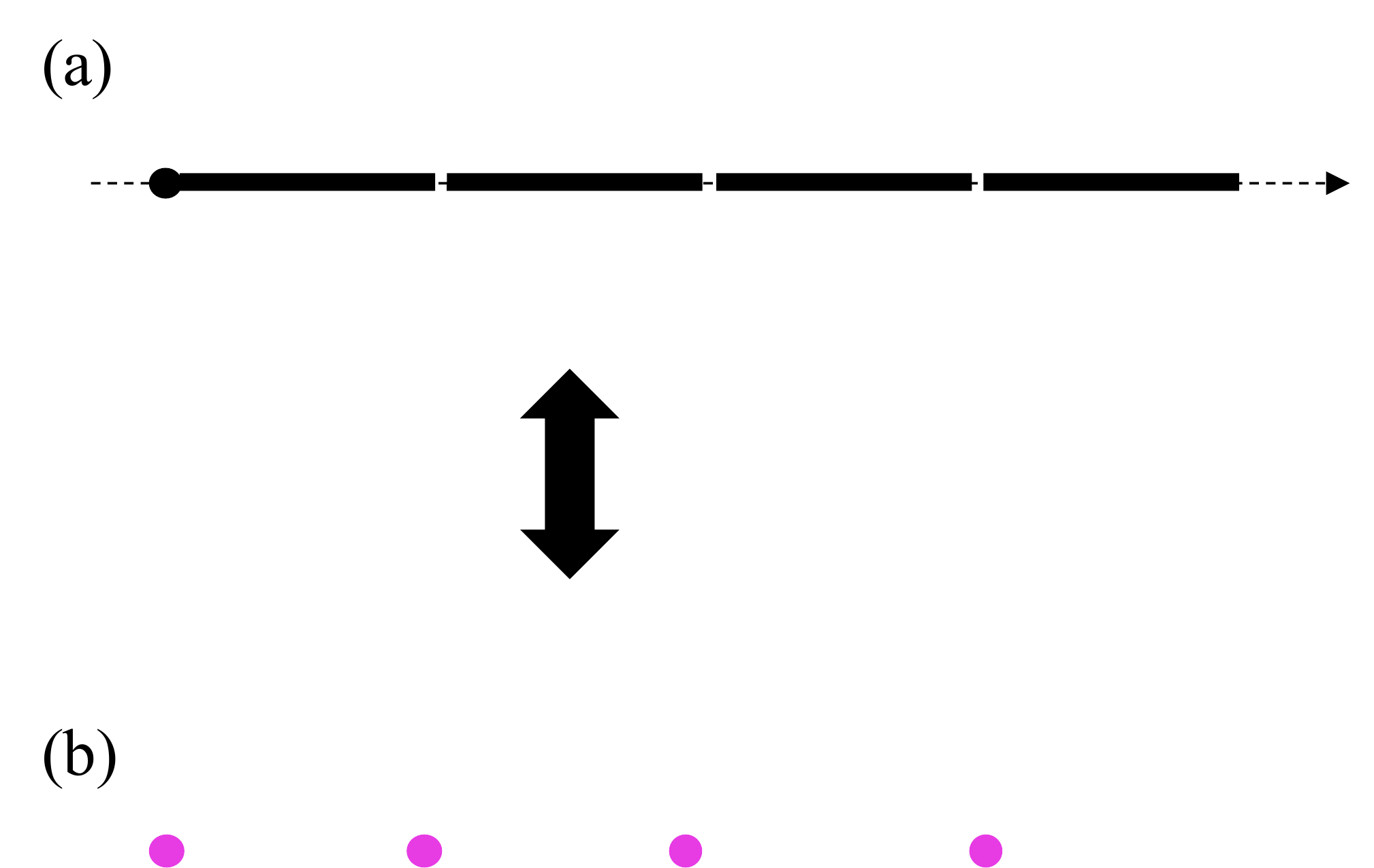}}}\hspace{5pt}
\caption{(a) Geometrical net of 2-cube. (b) Color representation of the geometrical net of 2-cube.}  \label{net_cube_2d}
\end{figure}

Also, the geometrical net of a unit square is represented by the four unit line segments in one-dimensional Euclidean space as shown in Fig.~\ref{net_cube_2d}(a). By using the pink coloring, the geometrical net of a unit square is represented by the four pink points in one-dimensional Euclidean space as shown in Fig.~\ref{net_cube_2d}(b).

The boundary of a square is constructed by (i) 4 edges and (ii) 4 vertices. The property (i) is checked by the fact that the colored net of a square is constructed by four colored points as shown in Fig.~\ref{net_cube_2d}(b)\footnote{We call this kind of figure a colored geometrical net.}. The property (ii) is checked by the theorem that Euler's characteristic is 0 for $2n+1$-dimensional regular polytopes.

\begin{figure}[H]
\centering
{%
\resizebox*{\textwidth}{!}{\includegraphics{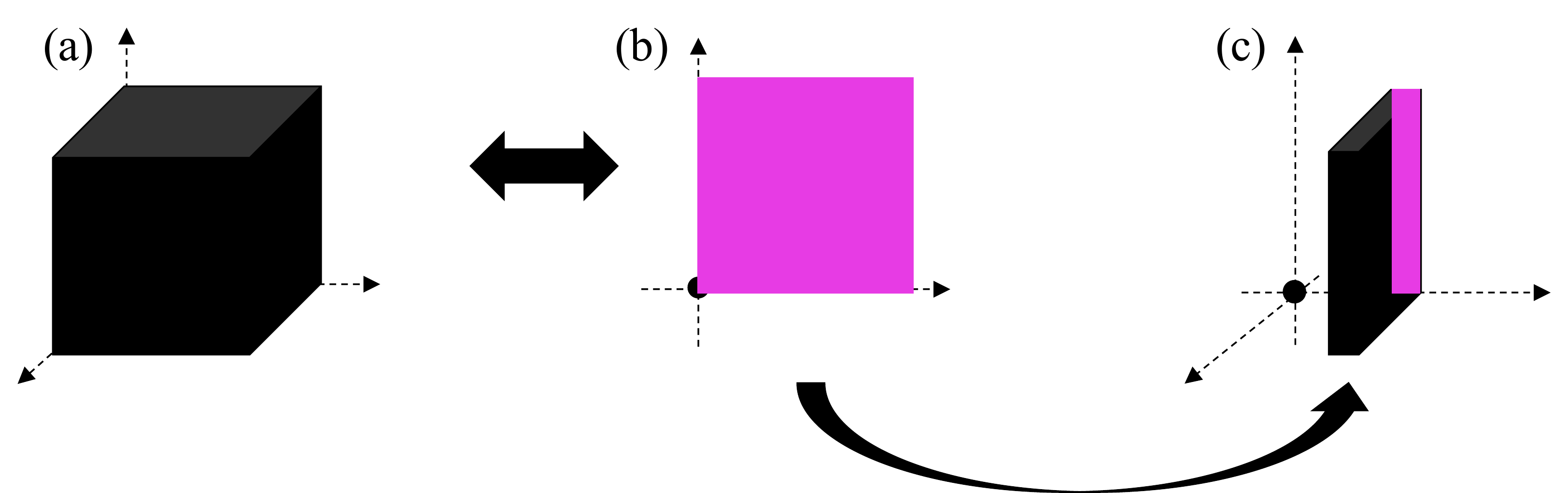}}}\hspace{5pt}
\caption{(a) 3-cube. (b) Color representation of 3-cube. (c) Decomposition of color representation of 3-cube.}  \label{cube_3d}
\end{figure}

A three-dimensional hypercube (3-cube) is a cube. A unit cube (Fig.~\ref{cube_3d}(a)) is represented by a pink unit square with length 1 on one edge (Fig.~\ref{cube_3d}(b)). This colored figure can be interpreted as an infinite set of pink unit line segments corresponding to the unit squares obtained by extending the pink line segments in the direction of the 3rd dimension as shown in Fig.~\ref{cube_3d}(c).

\begin{figure}[H]
\centering
{%
\resizebox*{\textwidth}{!}{\includegraphics{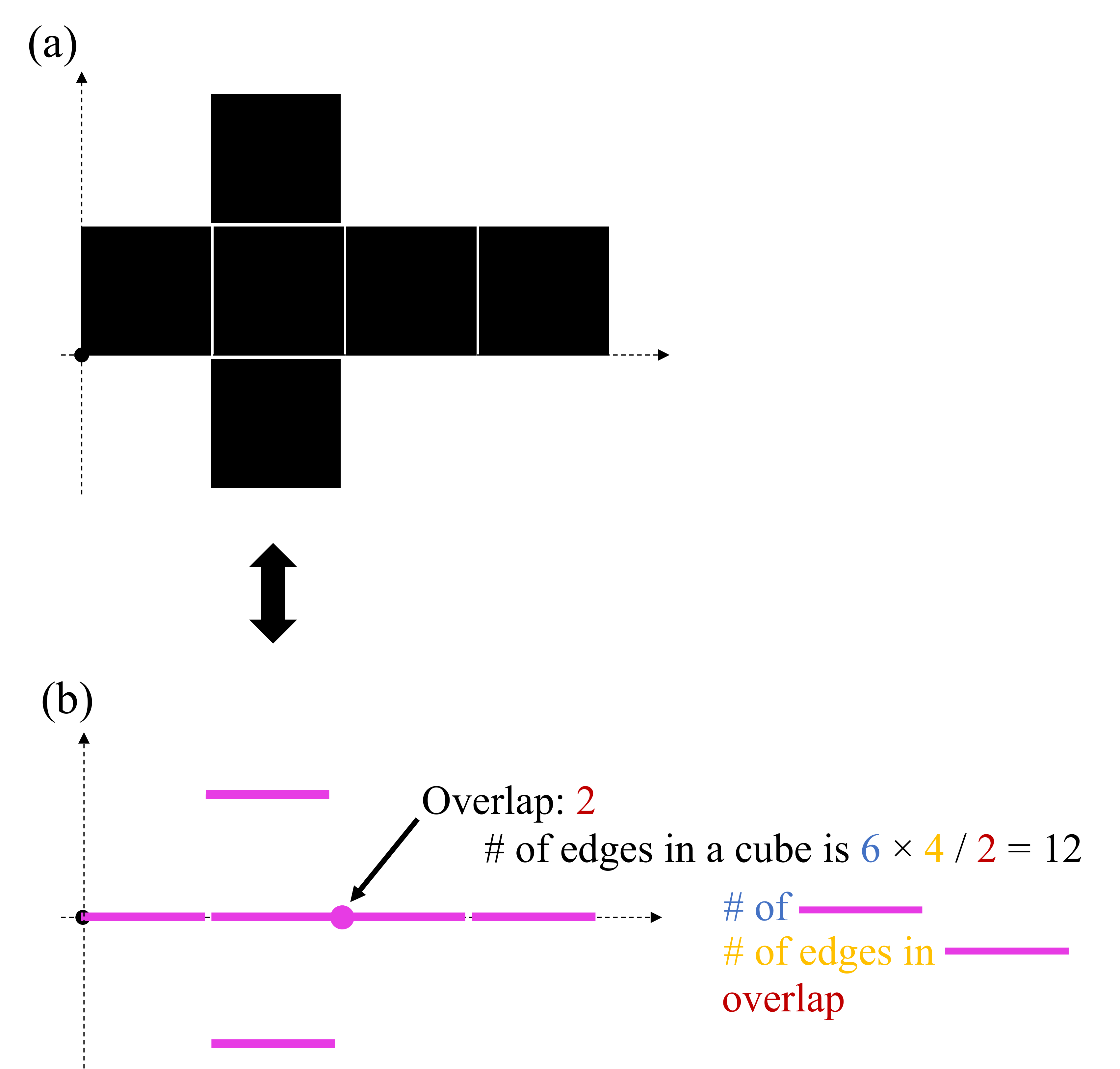}}}\hspace{5pt}
\caption{(a) Geometrical net of 3-cube. (b) Color representation of geometrical net of 3-cube.}  \label{net_cube_3d}
\end{figure}

Also, a geometrical net of a unit cube is represented by the six unit squares in two-dimensional Euclidean space as shown in Fig.~\ref{net_cube_3d}(a). By using the pink coloring, the geometrical net of a unit cube is represented by the eight pink line segments in two-dimensional Euclidean space as shown in Fig.~\ref{net_cube_3d}(b).

The boundary of a cube is constructed by (i) 6 faces, (ii) 12 edges, and (iii) 8 vertices. The property (i) is checked by the fact that the colored net of a cube is constructed by 6 colored line segments as shown in Fig.~\ref{net_cube_3d}(b). The property (ii) is checked by calculating 6 (the number of colored line segments in the colored net) times 4 (the number of edges of a square) divided by 2 (miximum number of overlap of vertices in the colored net) as shown in Fig.~\ref{net_cube_3d}(b). Note that this overlap corresponds to the overlap of the edges of the square in the boundary of a cube. We remark that one can discuss degree of elements of polytopes in $n$ dimensions in $n-2$ dimensions by using a colored geometrical net. The property (iii) is checked by the theorem that Euler's characteristic is 2 for $2n+1$-dimensional regular polytopes.

\begin{figure}[H]
\centering
{%
\resizebox*{\textwidth}{!}{\includegraphics{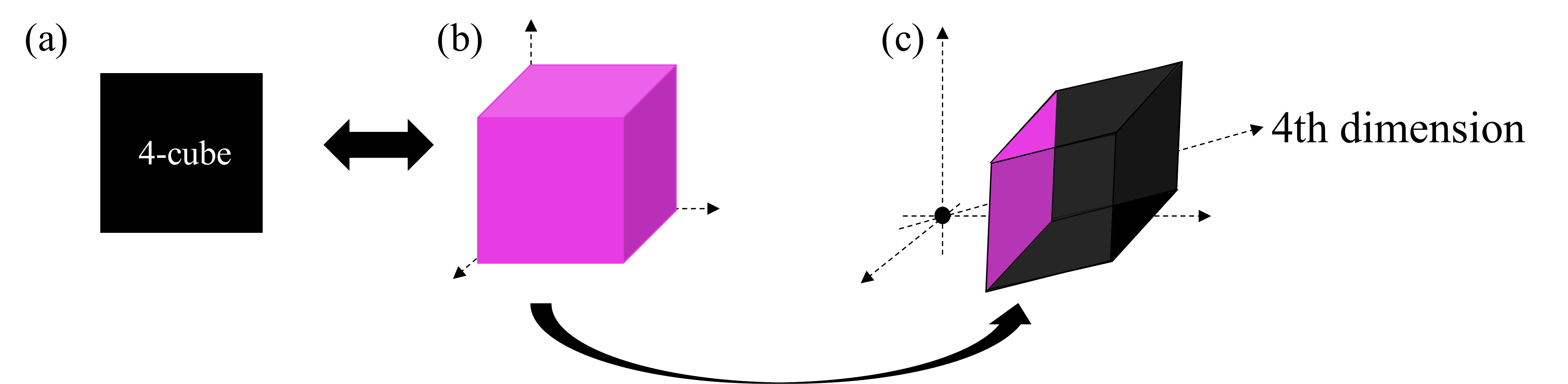}}}\hspace{5pt}
\caption{(a) 4-cube. (b) Color representation of 4-cube. (c) Decomposition of color representation of 4-cube.}  \label{cube_4d}
\end{figure}

A four-dimensional hypercube (4-cube) is a tesseract. A unit tesseract (Fig.~\ref{cube_4d}(a)) is represented by a pink unit cube with length 1 on one edge (Fig.~\ref{cube_4d}(b)). This colored figure can be interpreted as an infinite set of pink squares corresponding to the unit cubes obtained by extending the pink squares in the direction of the 4th dimension as shown in Fig.~\ref{cube_4d}(c).

\begin{figure}[H]
\centering
{%
\resizebox*{14cm}{!}{\includegraphics{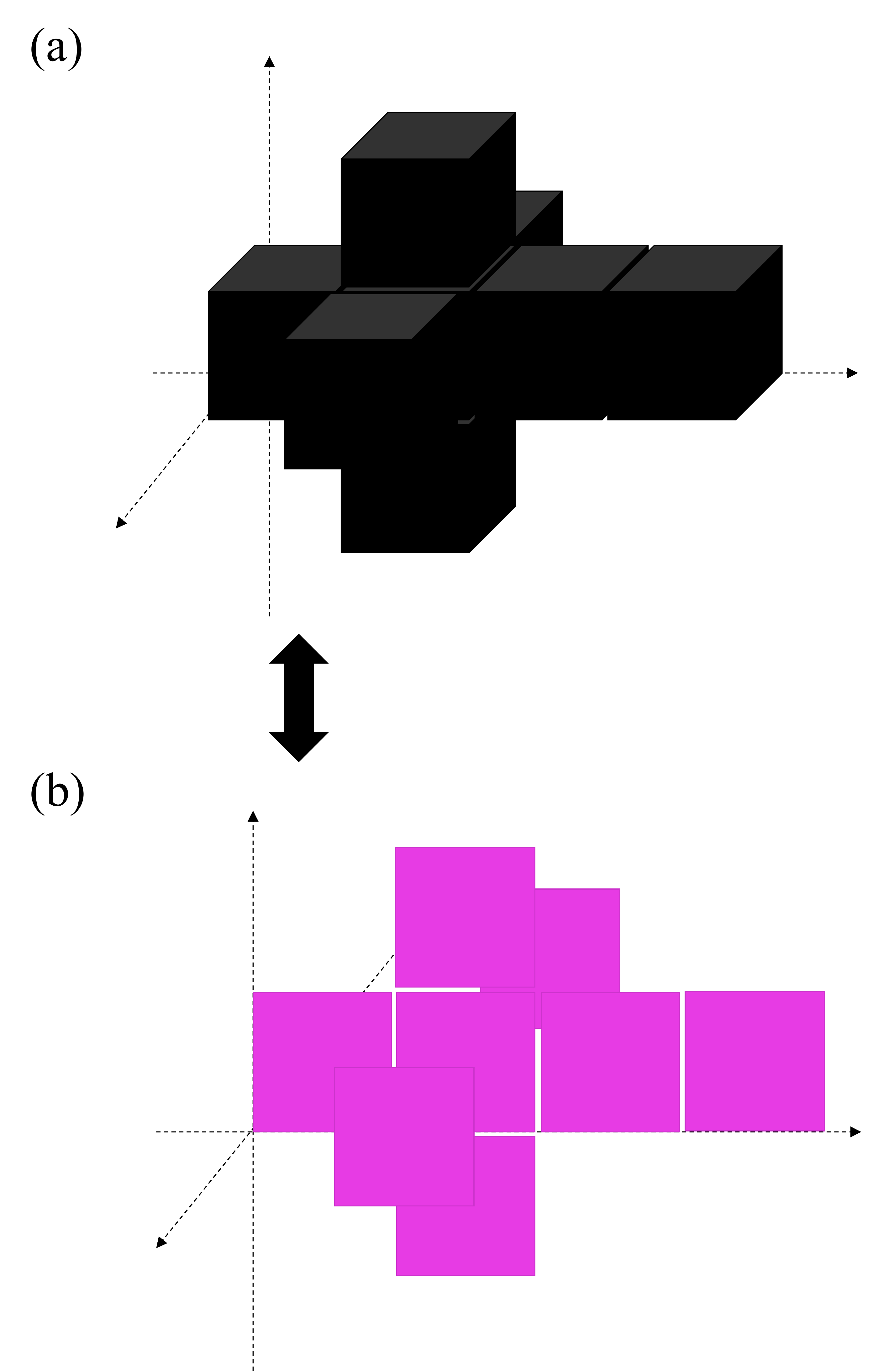}}}\hspace{5pt}
\caption{(a) Geometrical net of 4-cube. (b) Color representation of geometrical net of 4-cube.}  \label{net_cube_4d}
\end{figure}

Also, a geometrical net of a unit teseract is represented by the eight unit cubes in three-dimensional Euclidean space as shown in Fig.~\ref{net_cube_4d}(a). By using the pink coloring, the geometrical net of a unit tesseract is represented by the eight pink squares in three-dimensional Euclidean space as shown in Fig.~\ref{net_cube_4d}(b).

The boundary of a tesseract is constructed by (i) 8 cells, (ii) 24 faces, (iii) 32 edges, and (iv) 16 vertices. The property (i) is checked by the fact that the colored net of a tesseract is constructed by 8 colored squares as shown in Fig.~\ref{net_cube_4d}(b). The property (ii) is checked by calculating 8 (the number of colored squares in the colored net) times 6 (the number of faces of a cube) divided by 2 (miximum number of overlap of edges in the colored net) as shown in Fig.~\ref{net_cube_4d}(b). The property (iii) is checked by calculating 8 (the number of colored squares in the colored net) times 12 (the number of edges of a cube) divided by 3 (miximum number of overlap of vertices in the colored net) as shown in Fig.~\ref{net_cube_4d}(b). The property (iv) is checked by the theorem that Euler's characteristic is 0 for $2n+1$-dimensional regular polytopes.

\begin{figure}[H]
\centering
{%
\resizebox*{\textwidth}{!}{\includegraphics{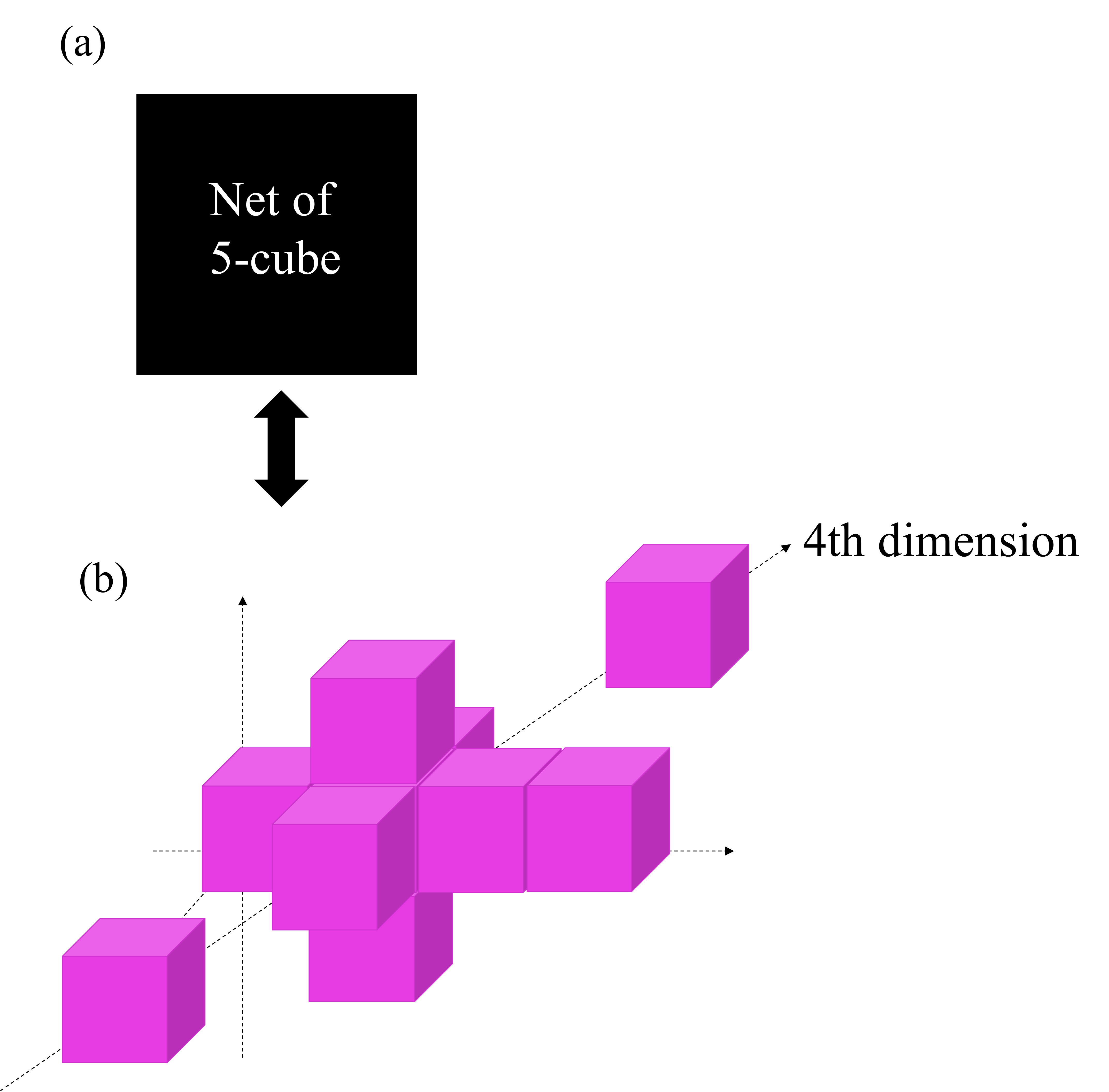}}}\hspace{5pt}
\caption{(a) Geometrical net of 5-cube. (b) Color representation of geometrical net of 5-cube.}  \label{net_cube_5d}
\end{figure}

Then, we propose a way to represent 5-cube by using the colored geometrical net. The geometrical net of a unit 5-cube is represented by the ten pink cubes in four-dimensional Euclidean space as shown in Fig.~\ref{net_cube_5d}.

The boundary of a 5-cube is constructed by (i) 10 4-faces, (ii) 40 cells, (iii) 80 faces, and (iv) 80 edges, and (v) 32 vertices. The property (i) is checked by the fact that the colored net of a 5-cube is constructed by 10 colored cube as shown in Fig.~\ref{net_cube_5d}. The property (ii) is checked by calculating 10 (the number of colored cubes in the colored net) times 8 (the number of cells of a tesseract) divided by 2 (miximum number of overlap of faces in the colored net) as shown in Fig.~\ref{net_cube_5d}. The property (iii) is checked by calculating 10 (the number of colored cubes in the colored net) times 24 (the number of faces of a tesseract) divided by 3 (miximum number of overlap of edges in the colored net) as shown in Fig.~\ref{net_cube_5d}. The property (iv) is checked by calculating 10 (the number of colored cubes in the colored net) times 32 (the number of edges of a tesseract) divided by 3 (miximum number of overlap of vertices in the colored net) as shown in Fig.~\ref{net_cube_5d}. The property (v) is checked by the theorem that Euler's characteristic is 2 for $2n+1$-dimensional regular polytopes.

\subsection{Hypertetrahedron}\label{Hypertetrahedron}
A one-dimensional hypertetrahedron (1-simplex) is a line segment.

\begin{figure}[H]
\centering
{%
\resizebox*{\textwidth}{!}{\includegraphics{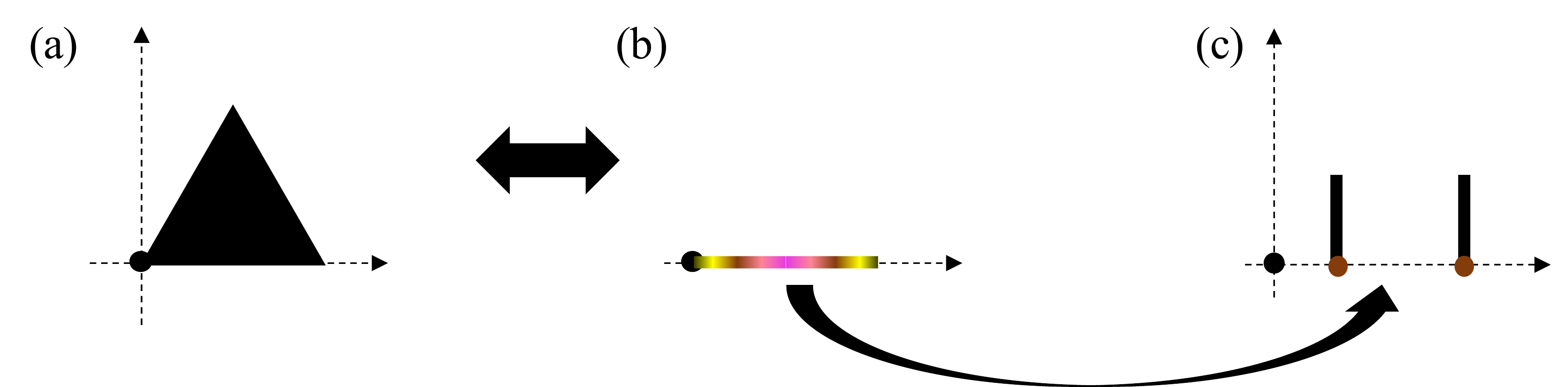}}}\hspace{5pt}
\caption{(a) 2-simplex. (b) Color representation of 2-simplex. (c) Decomposition of color representation of 2-simplex.}  \label{simplex_2d}
\end{figure}

A two-dimensional hypertetrahedron (2-simplex) is an equilateral triangle. A unit equilateral triangle (Fig.~\ref{simplex_2d}(a)) is represented by a black-pink-black linearly-graduated line segment whose length is $1$ (Fig.~\ref{simplex_2d}(b))\footnote{We call this gradated coloring simplex coloring. When using the simplex coloring, we set the maximum value in the color bar to $\sqrt{3}/2$, not 1.}. This colored figure can be interpreted as an infinite level sets of line segment boundaries as shown in Fig.~\ref{simplex_2d}(c).

\begin{figure}[H]
\centering
{%
\resizebox*{8cm}{!}{\includegraphics{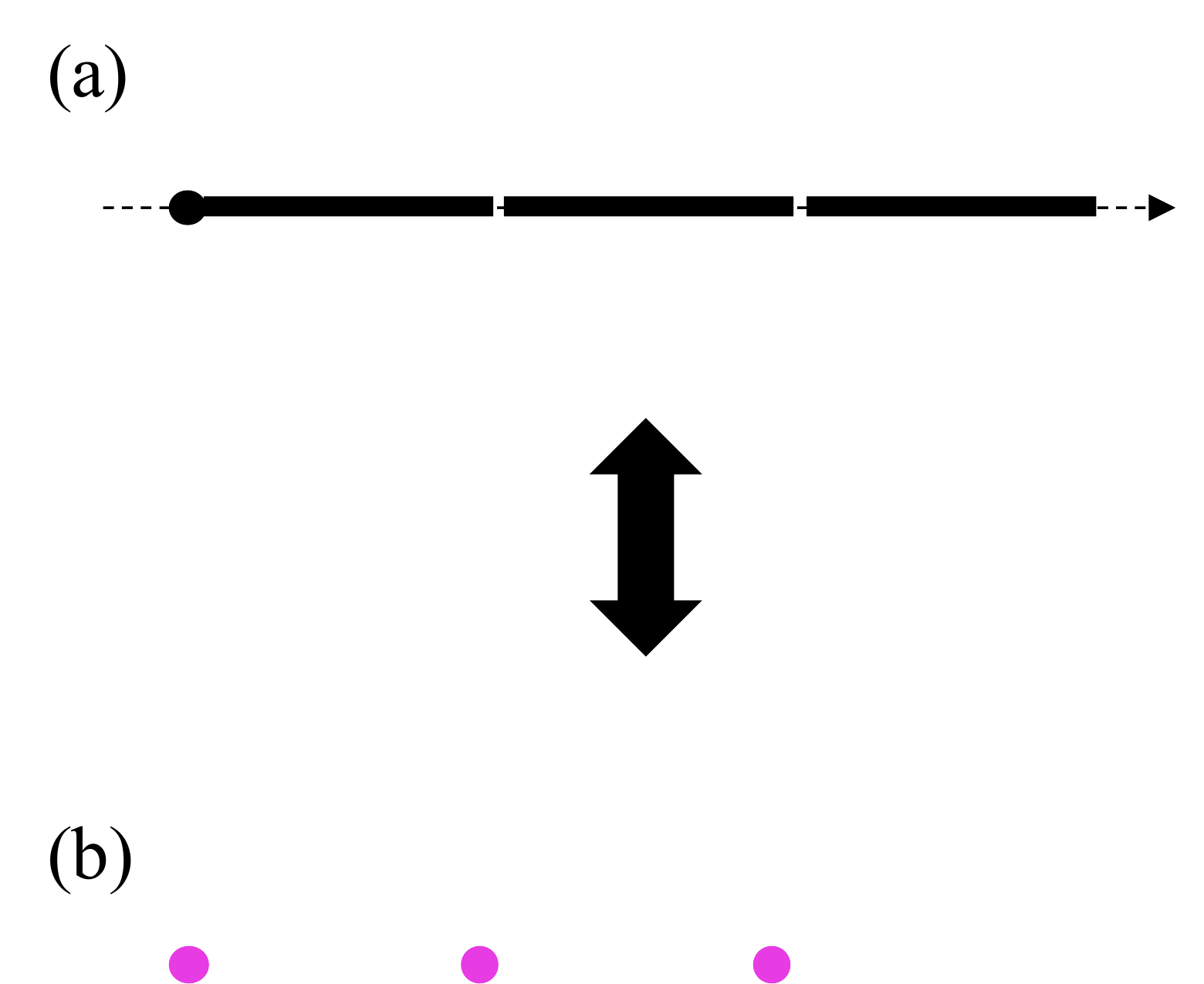}}}\hspace{5pt}
\caption{(a) Geometrical net of 2-simplex. (b) Color representation of geometrical net of 2-simplex.}  \label{net_simplex_2d}
\end{figure}

Also, a geometrical net of a unit equilateral triangle is represented by the three unit line segments in one-dimensional Euclidean space as shown in Fig.~\ref{net_simplex_2d}(a). By using the black-pink-black gradated coloring and its reverse coloring (black-green-black coloring), the geometrical net of a unit equilateral triangle is represented by the one black-pink-black gradated point and the two black-pink-black gradated points in one-dimensional Euclidean space with folding as shown in Fig.~\ref{net_simplex_2d}(b).

The boundary of an equilateral triangle is constructed by (i) 3 edges and (ii) 3 vertices. The property (i) is checked by the fact that the colored net of an equilateral triangle is constructed by 3 colored points as shown in Fig.~\ref{net_simplex_2d}. The property (ii) is checked by the theorem that Euler's characteristic is 0 for $2n+1$-dimensional regular polytopes.

\begin{figure}[H]
\centering
{%
\resizebox*{\textwidth}{!}{\includegraphics{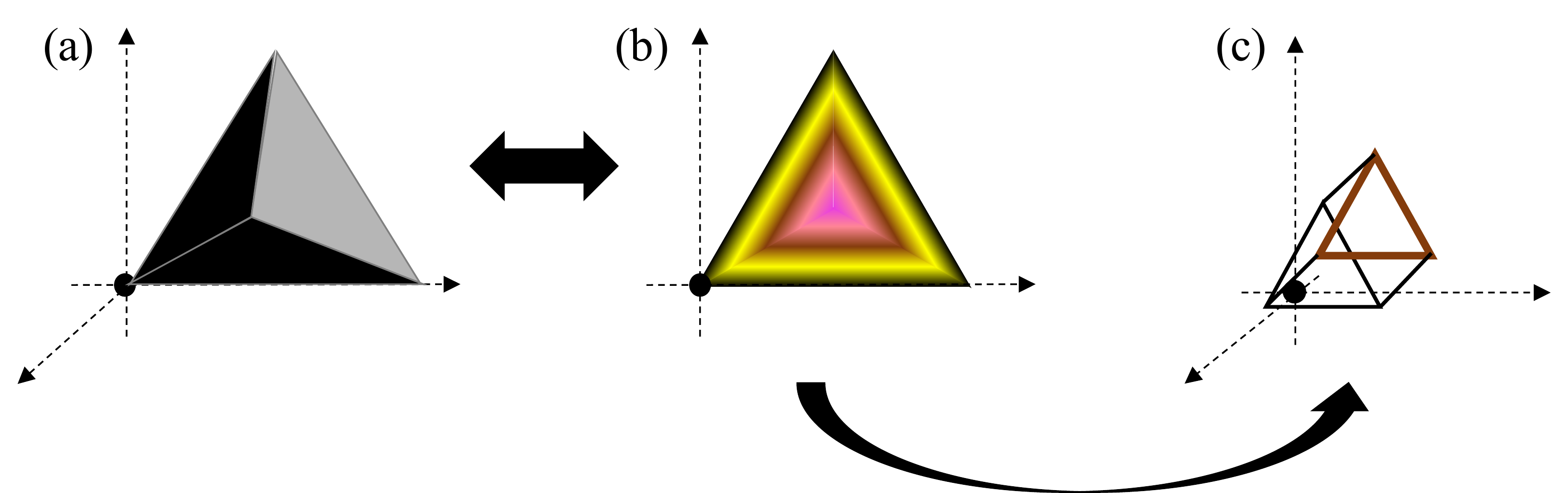}}}\hspace{5pt}
\caption{(a) 3-simplex. (b) Color representation of 3-simplex. (c) Decomposition of color representation of 3-simplex.}  \label{simplex_3d}
\end{figure}

A three-dimensional hypertetrahedron (3-simplex) is a tetrahedron. A unit regular tetrahedron (Fig.~\ref{simplex_3d}(a)) is represented by a black-pink-black gradated unit equilateral triangle with length 1 on one edge (Fig.~\ref{simplex_3d}(b)). This colored figure can be interpreted as an infinite level set of triangle boundaries as shown in Fig.~\ref{simplex_3d}(c). 

\begin{figure}[H]
\centering
{%
\resizebox*{8cm}{!}{\includegraphics{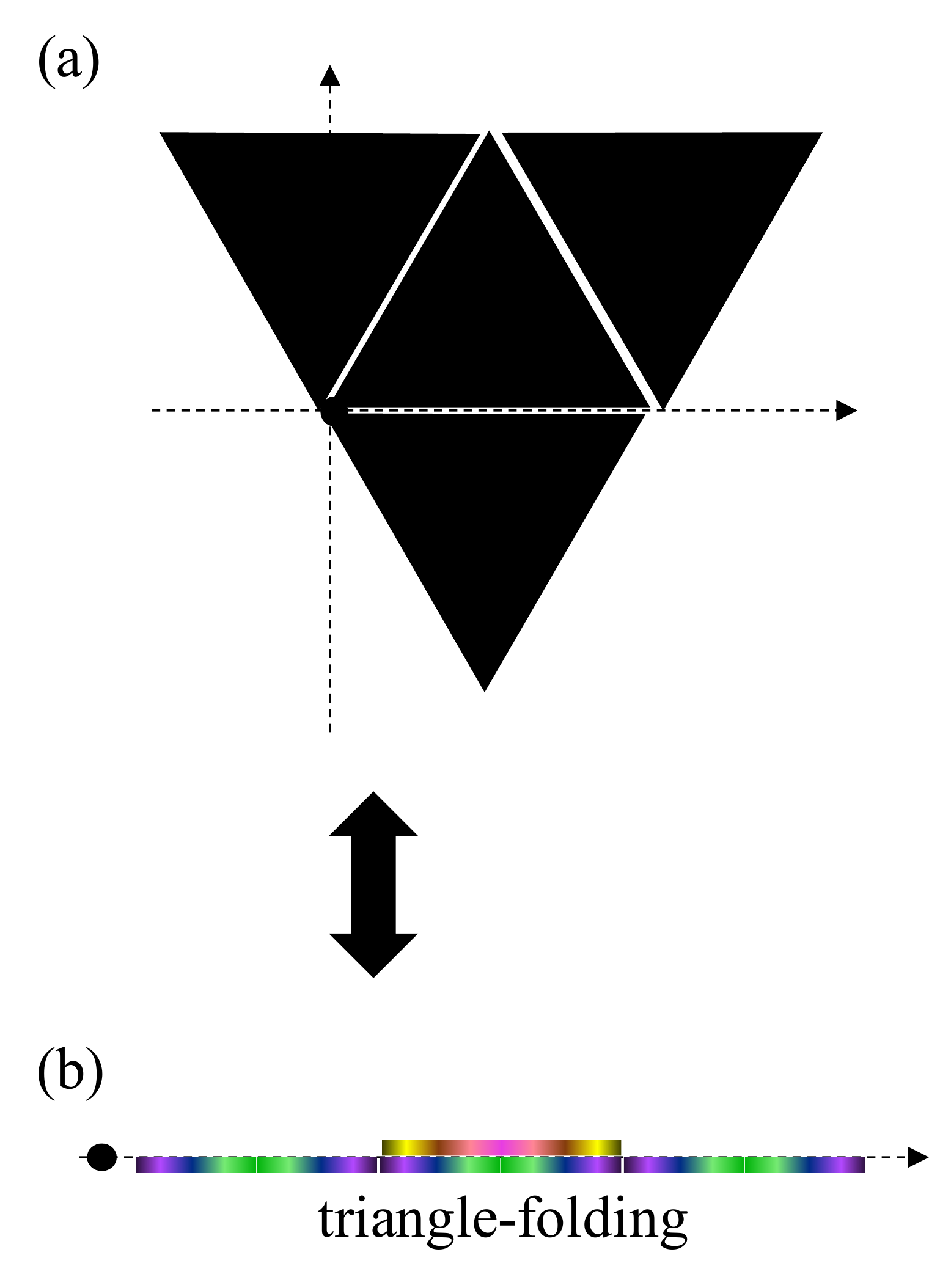}}}\hspace{5pt}
\caption{(a) Geometrical net of 3-simplex. (b) Color representation of geometrical net of 3-simplex.}  \label{net_simplex_3d}
\end{figure}

Also, a geometrical net of a regular tetrahedron is represented by the four unit equilateral triangles in two-dimensional Euclidean space as shown in Fig.~\ref{net_simplex_3d}(a). By using the black-pink/green-black gradated coloring, the geometrical net of a regular tetrhedron is represented by the one black-pink-black gradated line segment and the black-green-black gradated line segments in two-dimensional Euclidean space as shown in Fig.~\ref{net_simplex_3d}(b).

The boundary of a regular tetrahedron is constructed by (i) 4 faces, (ii) 6 edges, and (iii) 4 vertices. The property (i) is checked by the fact that the colored net of a tetrahedron is constructed by 4 colored line segments as shown in Fig.~\ref{net_simplex_3d}. The property (ii) is checked by calculating 4 (the number of black-pink-black gradated line segments in the colored net) times 3 (the number of edges of an equilateral triangle) divided by 2 (miximum number of overlap of vertices in the colored net) as shown in Fig.~\ref{net_simplex_3d}. The property (iii) is checked by the theorem that Euler's characteristic is 2 for $2n+1$-dimensional regular polytopes.

\begin{figure}[H]
\centering
{%
\resizebox*{\textwidth}{!}{\includegraphics{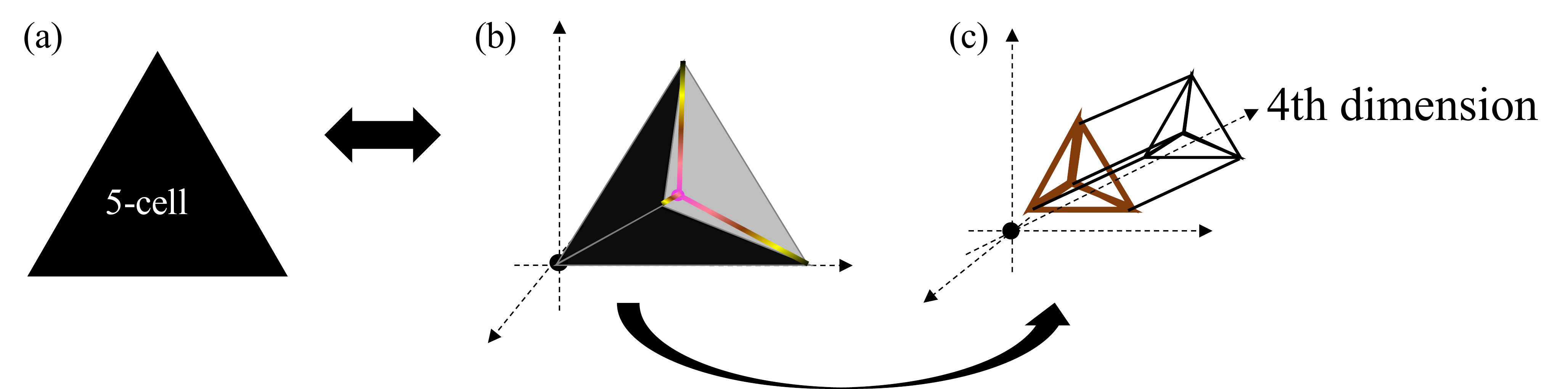}}}\hspace{5pt}
\caption{(a) 4-simplex. (b) Color representation of 4-simplex. (c) Decomposition of color representation of 4-simplex.}  \label{simplex_4d}
\end{figure}

A four-dimensional hypertetrahedron (4-simplex) is a 5-cell. A unit 5-cell (Fig.~\ref{simplex_4d}(a)) is represented by a black-pink-black gradated regular tetrahedron with length 1 on one edge (Fig.~\ref{simplex_4d}(b)). Note that colored lines in Fig.~\ref{simplex_4d}(b) is auxiliary lines for us living in three-dimensional world. This colored figure can be interpreted as an infinite set of tetrahedron contours as shown in Fig.~\ref{simplex_4d}(c).

\begin{figure}[H]
\centering
{%
\resizebox*{\textwidth}{!}{\includegraphics{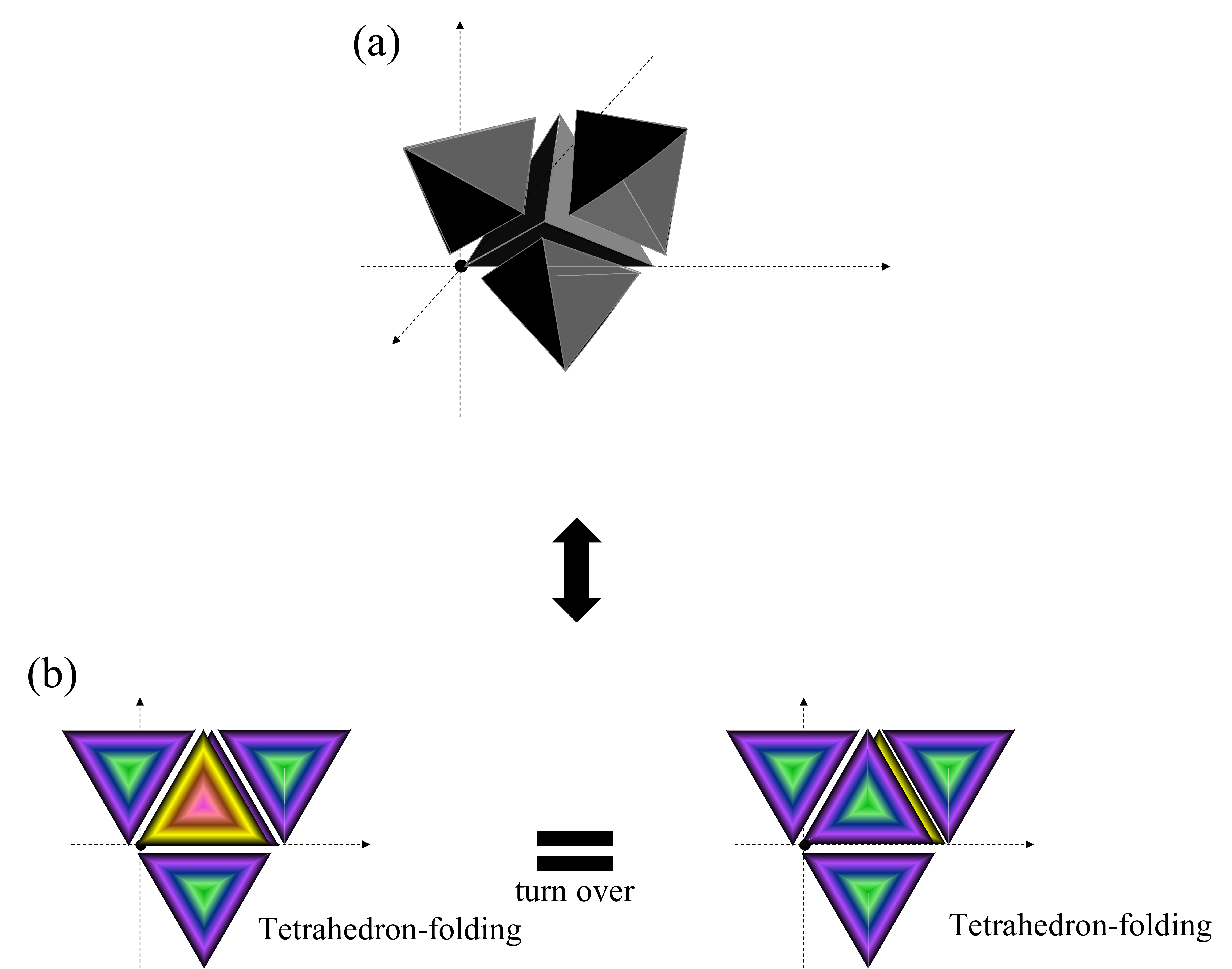}}}\hspace{5pt}
\caption{(a) Geometrical net of 4-simplex. (b) Color representation of geometrical net of 4-simplex where the center triangle is double.}  \label{net_simplex_4d}
\end{figure}

Also, a geometrical net of a unit 5-cell is represented by the five regular tetrahedra in three-dimensional Euclidean space as shown in Fig.~\ref{net_simplex_4d}(a). By using the black-pink/green-black gradated coloring, the geometrical net of a unit 5-cell is represented by the one black-pink-black gradated triangle and the four black-green-black gradated triangles in three-dimensional Euclidean space as shown in Fig.~\ref{net_simplex_4d}(b).

The boundary of a 5-cell is constructed by (i) 5 cells, (ii) 10 faces, (iii) 10 edges, and (iv) 5 vertices. The property (i) is checked by the fact that the colored net of a 5-cell is constructed by 5 colored triangles in Fig.~\ref{net_simplex_4d}. The property (ii) is checked by calculating 5 (the number of colored triangles in the colored net) times 4 (the number of faces of a tetrahedron) divided by 2 (miximum number of overlap of edges in the colored net) as shown in Fig.~\ref{net_simplex_4d}. The property (iii) is checked by calculating 5 (the number of colored squares in the colored net) times 6 (the number of edges of a tetrahedron) divided by 3 (miximum number of overlap of vertices in the colored net) as shown in Fig.~\ref{net_simplex_4d}. The property (iv) is checked by the theorem that Euler's characteristic is 0 for $2n+1$-dimensional regular polytopes.

\begin{figure}[H]
\centering
{%
\resizebox*{\textwidth}{!}{\includegraphics{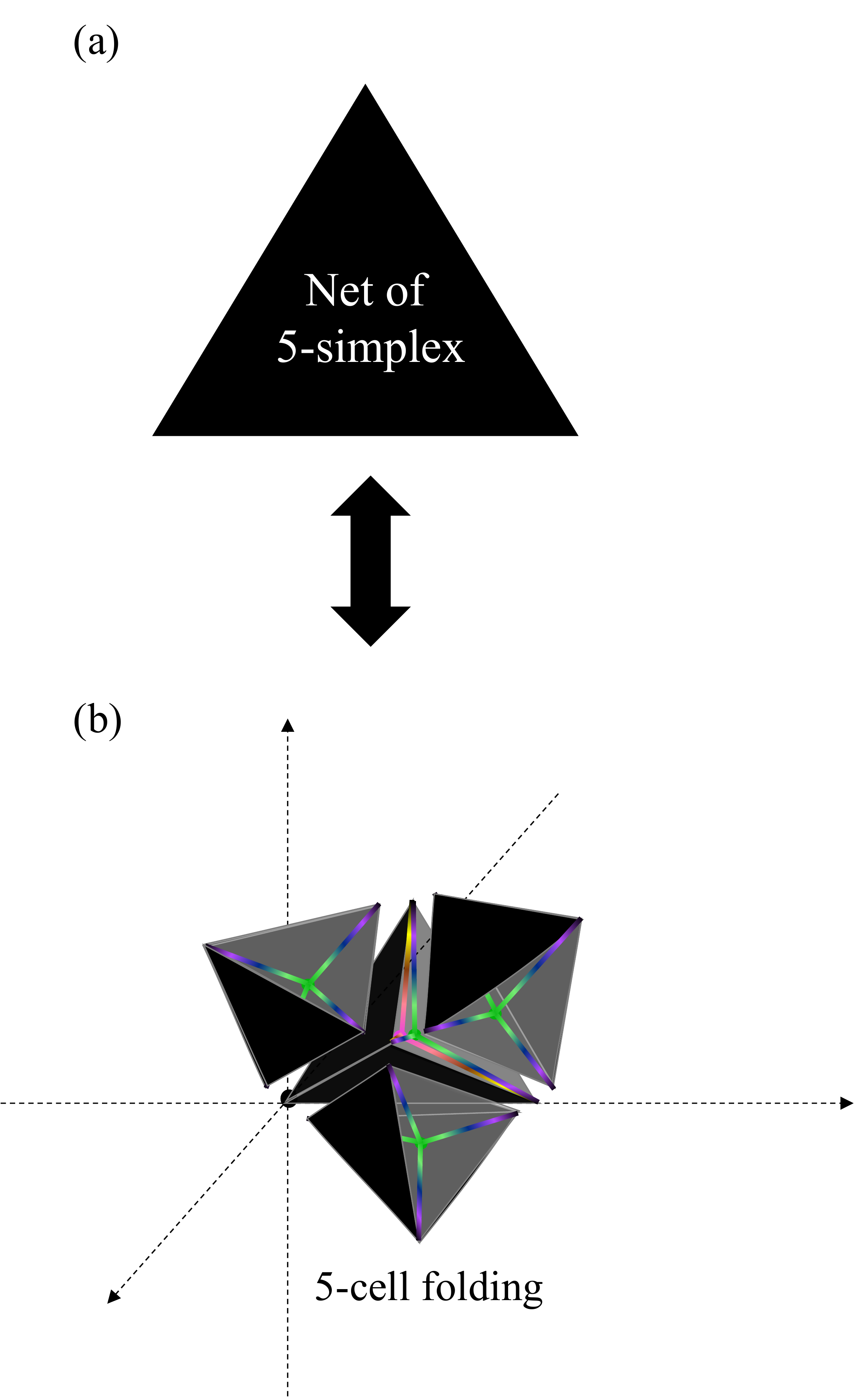}}}\hspace{5pt}
\caption{(a) Geometrical net of 5-simplex. (b) Color representation of geometrical net of 5-simplex where the center tetrahedron is double.}  \label{net_simplex_5d}
\end{figure}

Then, we propose a way to represent 5-simplex by using the colored geometrical net. The geometrical net of a unit 5-simplex is represented by the one black-pink-black gradated regular tetrahedron the five black-geen-black gradated regular tetrahedra in four-dimensional Euclidean space as shown in Fig.~\ref{net_simplex_5d}. 

The boundary of a 5-simplex is constructed by (i) 6 4-faces, (ii) 15 cells, (iii) 20 faces, and (iv) 15 edges, and (v) 6 vertices. The property (i) is checked by the fact that the colored net of a 5-simplex is constructed by 6 colored tetrahedra as shown in Fig.~\ref{net_simplex_5d}. The property (ii) is checked by calculating 6 (the number of colored tetrahedra in the colored net) times 5 (the number of cells of a 5-cell) divided by 2 (miximum number of overlap of faces in the colored net) as shown in Fig.~\ref{net_simplex_5d}. The property (iii) is checked by calculating 6 (the number of colored tetrahedra in the colored net) times 10 (the number of faces of a 5-cell) divided by 3 (miximum number of overlap of edges in the colored net) as shown in Fig.~\ref{net_simplex_5d}. The property (iv) is checked by calculating 6 (the number of colored tetrahedra in the colored net) times 10 (the number of edges of a 5-cell) divided by 4 (miximum number of overlap of vertices in the colored net) as shown in Fig.~\ref{net_simplex_5d}. The property (v) is checked by the theorem that Euler's characteristic is 2 for $2n+1$-dimensional regular polytopes.

\newpage
\section{Hyperprism}\label{Hyperprism}
In this section, we propose a way to represent some hyperprisms by using the cube coloring and/or the simplex coloring.

\subsection{Polytopal prism}
\begin{figure}[H]
\centering
{%
\resizebox*{6cm}{!}{\includegraphics{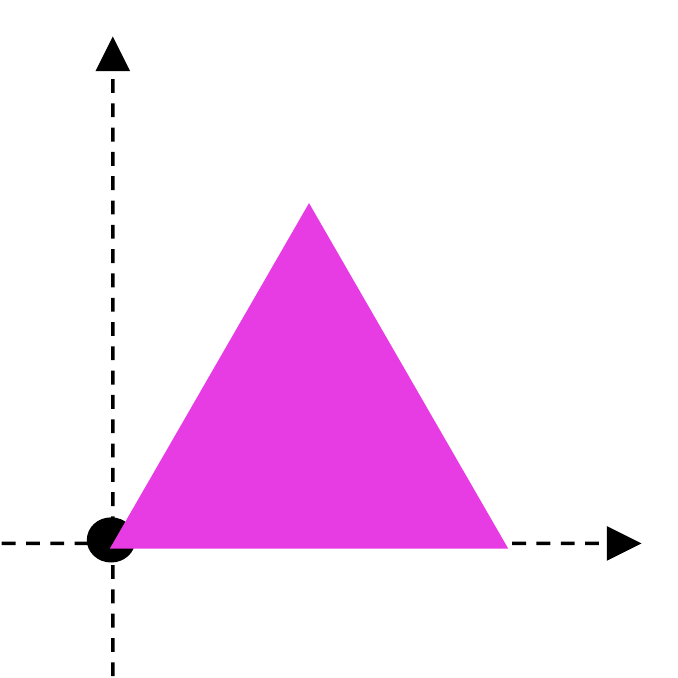}}}\hspace{5pt}
\caption{Color representation of a polytopal prism constructed by the Cartesian product of a 2-simplex and a unit line segment.}\label{prism_2-simplex}
\end{figure}

\begin{figure}[H]
\centering
{%
\resizebox*{8cm}{!}{\includegraphics{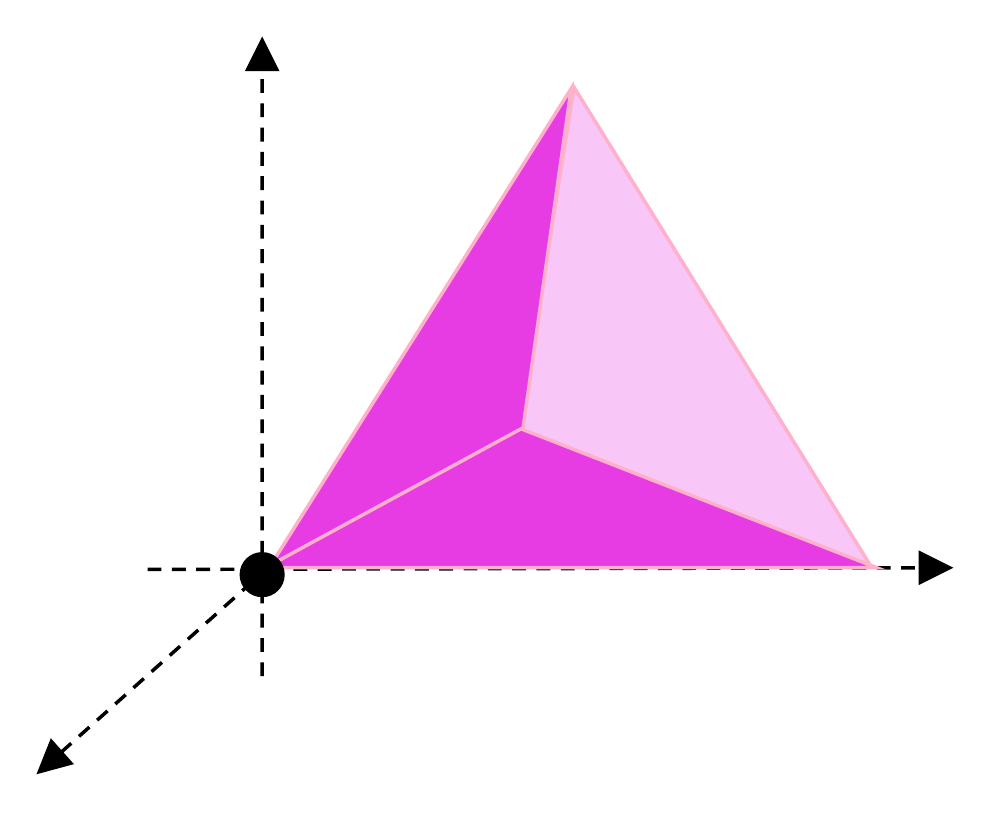}}}\hspace{5pt}
\caption{Color representation of a polytopal prism constructed by the Cartesian product of a 3-simplex and a unit line segment.}\label{prism_3-simplex}
\end{figure}

We define a polytopal prism as the Cartesian product of an $n$-dimensional polytope ($n\geq2$) and a unit line segment. Because a unit line segment can be represented by a pink point, then, a polytopal prism constructed by $n$-dimensional polytope $X$ and a unit line segment can be represented by the $n$-dimensional polytope colored by pink. The color representation of $n$-cubes shown in Sec.~\ref{Hypercube} are the simplest examples of the coloring of polytopal prisms.
In addition, one can construct a polytopal prism by the Cartesian product of an $n$-simplex and a unit line segment. Fig.~\ref{prism_2-simplex} shows the color representation of a polytopal prism constructed by the Cartesian product of a 2-simplex (equilateral triangle) and a unit line segment, i.e., triangular prism. Fig.~\ref{prism_3-simplex} shows the color representation of a polytopal prism constructed by the Cartesian product of a 3-simplex (tetrahedron) and a unit line segment, i.e., tetrahedral prism.

\subsection{Duoprism}
In the previous subsection, we dealt with polytopal prisms that are the Cartesian product of an $n$-dimensional polytope ($n\geq2$) and a line segment. In this subsection, we deal with duoprisms that are the Cartesian product of two $n$-dimensional polytopes ($n\geq2$). Especially, we study $p$-$q$ duoprisms constructed by a $p$-gon and $q$-gon in four dimensions.

\begin{figure}[H]
\centering
{%
\resizebox*{6cm}{!}{\includegraphics{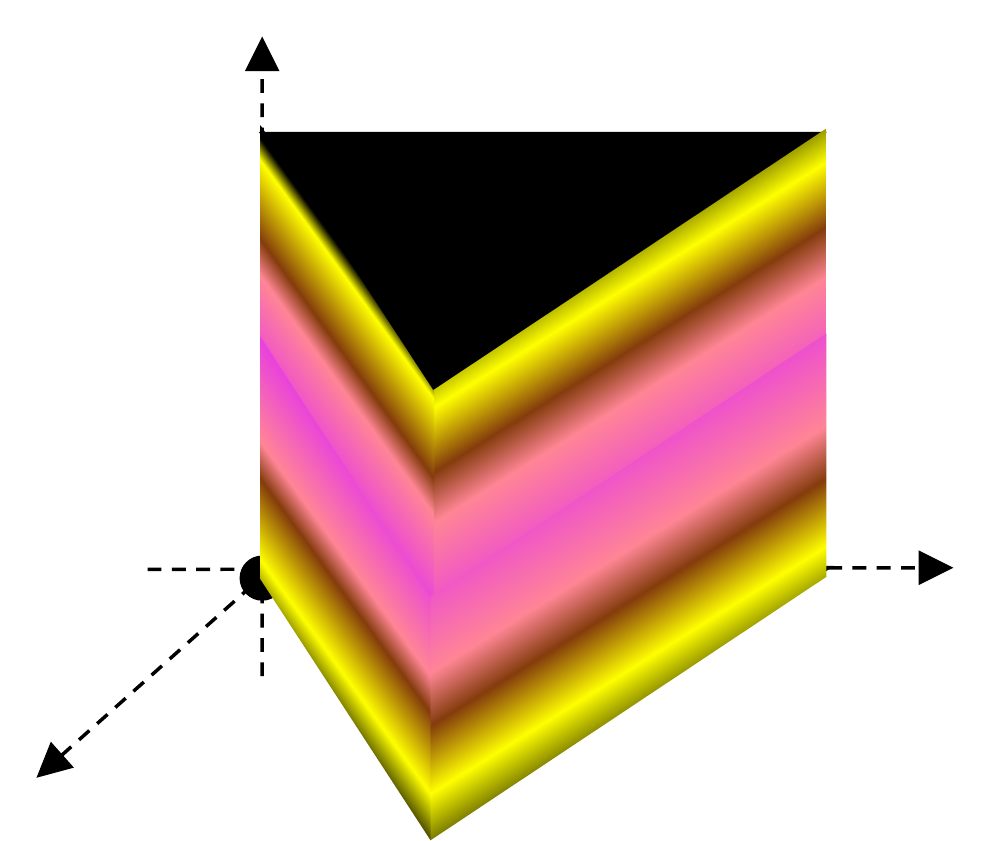}}}\hspace{5pt}
\caption{Color representation of 3-3 duoprism.}\label{prism_3-3}
\end{figure}

\begin{figure}[H]
\centering
{%
\resizebox*{6cm}{!}{\includegraphics{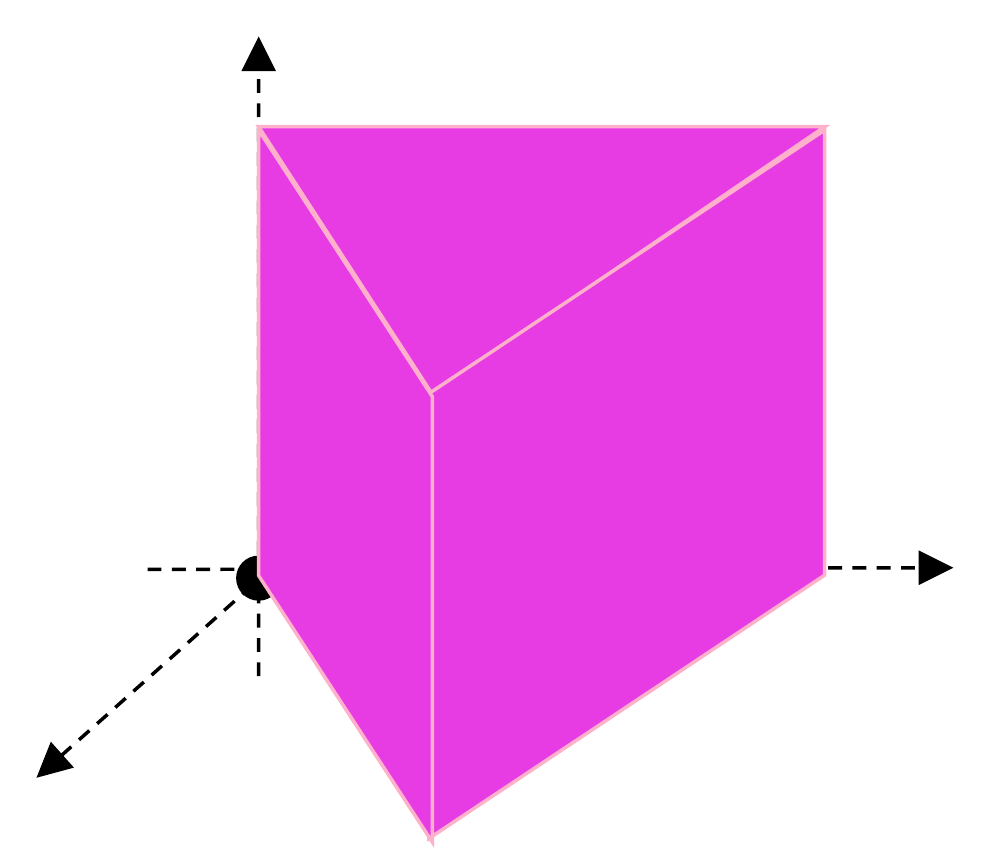}}}\hspace{5pt}
\caption{Color representation of 3-4 duoprism.}\label{prism_3-4}
\end{figure}

\begin{figure}[H]
\centering
{%
\resizebox*{6cm}{!}{\includegraphics{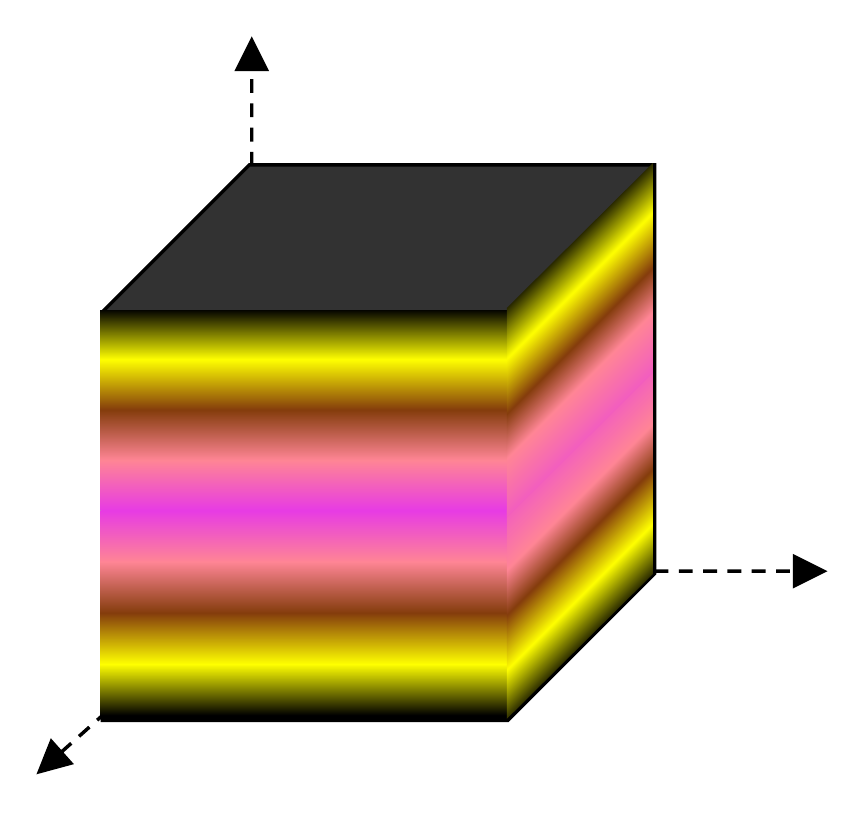}}}\hspace{5pt}
\caption{Color representation of 4-3 duoprism.}\label{prism_4-3}
\end{figure}

Fig.~\ref{prism_3-3} shows the color representation of a 3-3 duoprism. We remark that this figure is constructed by an equilateral triangle and the color representation of an equilateral triangle that can be represented by a black-pink-black gradated line segment, and then, one can represent a 3-3 duoprism that is a four-dimensional geometrical object in three-dimensional Euclidean space. Similarly, one can obtain the color representation of a 3-4 duoprism (Fig.~\ref{prism_3-4}).  We remark that this figure is constructed by an equilateral triangle and the color representation of a square that can be represented by a pink line segment. This colored figure can also be regarded as a polytopal prism constructed by a triangular prism and a line segment (pink point). Also, one can obtain the color representation of a 4-3 duoprism (Fig.~\ref{prism_4-3}) that is equivalent to a 3-4 duoprism. We remark that this figure is constructed by an square and the color representation of an equilateral triangle that can be represented by a black-pink-black gradated line segment. Fig.~\ref{cube_4d} shows the color representation of a 4-4 duoprism, i.e., a tesseract. We remark that this figure is constructed by a square and the color representation of a square that can be represented by a pink line segment.

\newpage
\section{Uncoloring and Truncation}\label{Uncoloring and Truncation}
In this section, we propose a way to represent a truncated polytope by using the uncoloring introduced in Sec.~\ref{Introduction}.

\begin{figure}[H]
\centering
{%
\resizebox*{\textwidth}{!}{\includegraphics{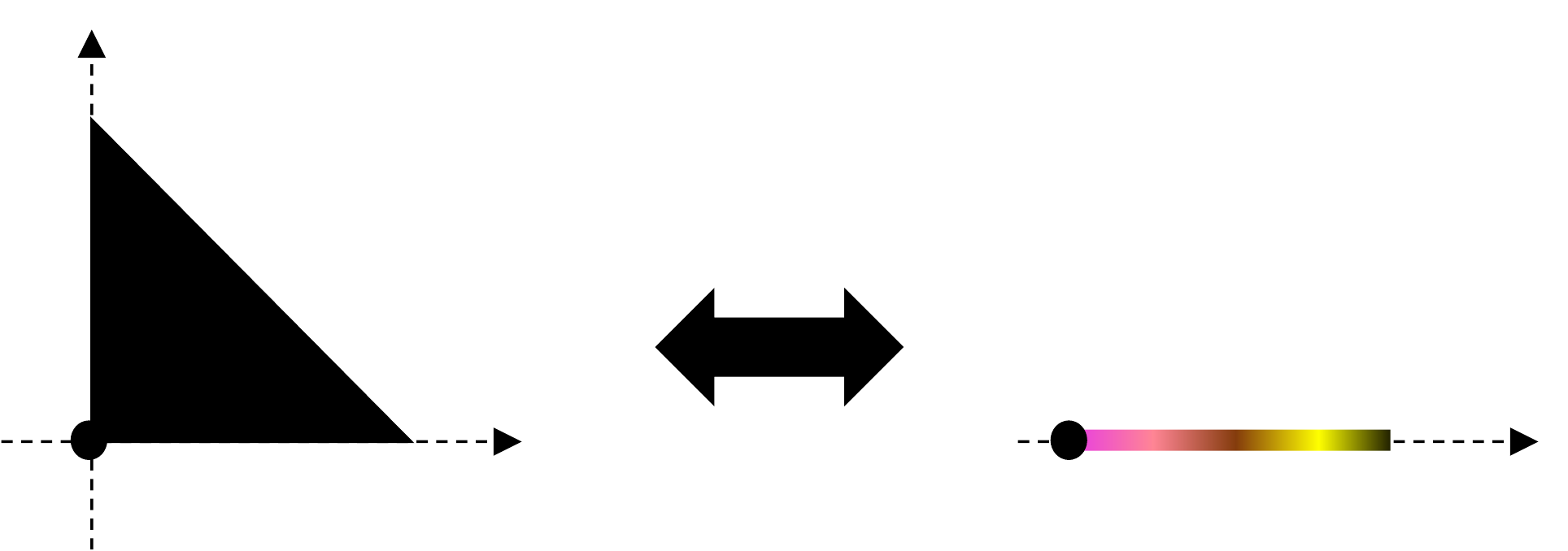}}}\hspace{5pt}
\caption{Color representation of 2-cube corner.}\label{corner_2d}
\end{figure}

\begin{figure}[H]
\centering
{%
\resizebox*{\textwidth}{!}{\includegraphics{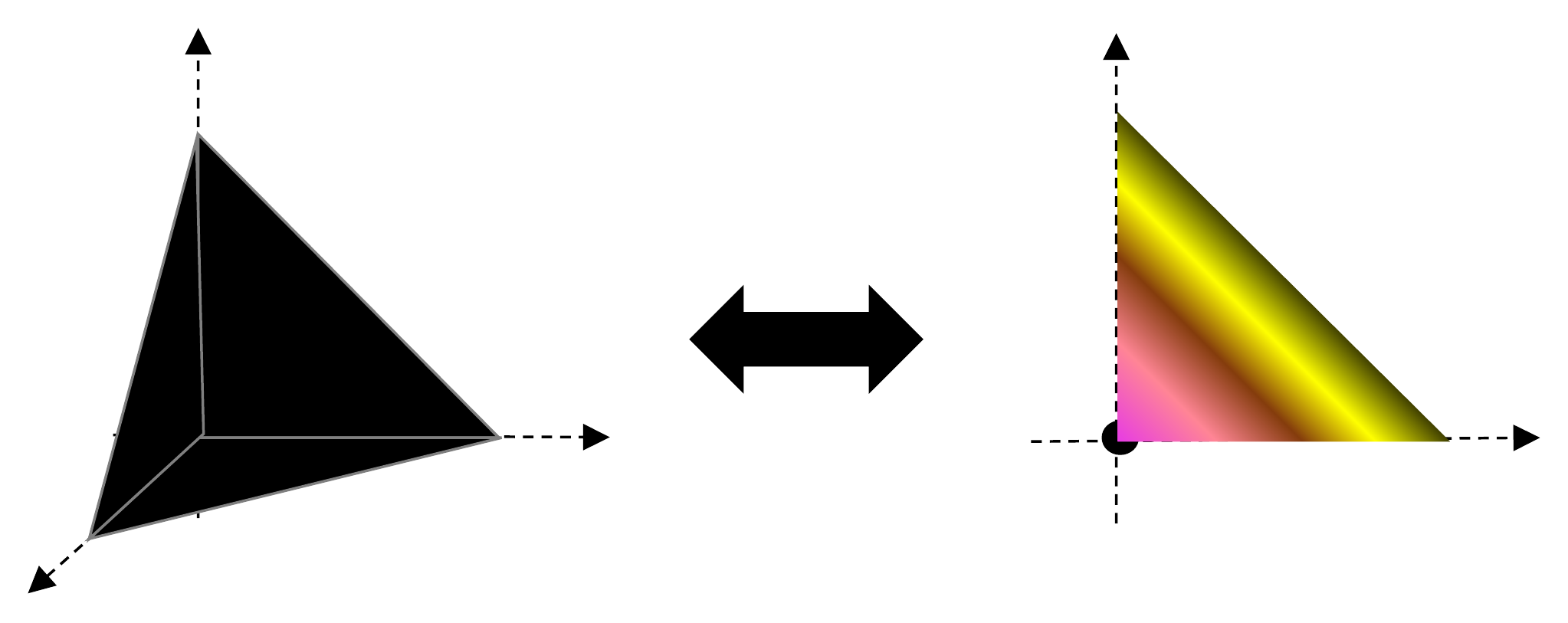}}}\hspace{5pt}
\caption{Color representation of 3-cube corner.}\label{corner_3d}
\end{figure}

\begin{figure}[H]
\centering
{%
\resizebox*{\textwidth}{!}{\includegraphics{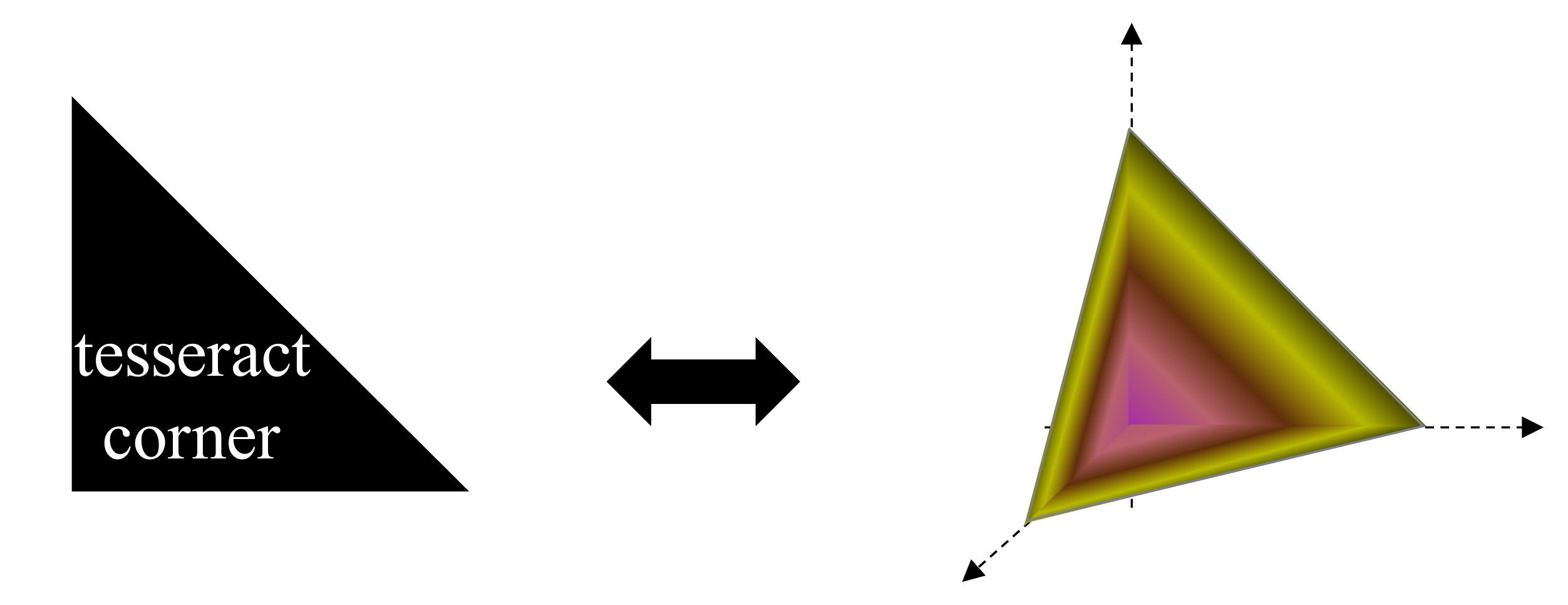}}}\hspace{5pt}
\caption{Color representation of 4-cube corner.}\label{corner_4d}
\end{figure}

Before proposing the color representation of truncated polytopes, we propose the color representation of the $n$-cube corner ($n\geq2$). Fig.~\ref{corner_2d} shows the color representation of the 2-cube corner (isosceles right triangle). One can obtain the color representation of the 3-cube corner (trirectangular tetrahedron) and find that it is represented by the coloring of the 2-cube corner (Fig.~\ref{corner_3d}). Similarly, the color representation of the 4-cube corner (Fig.~\ref{corner_4d}) is given by the coloring of the 3-cube corner.

\begin{figure}[H]
\centering
{%
\resizebox*{\textwidth}{!}{\includegraphics{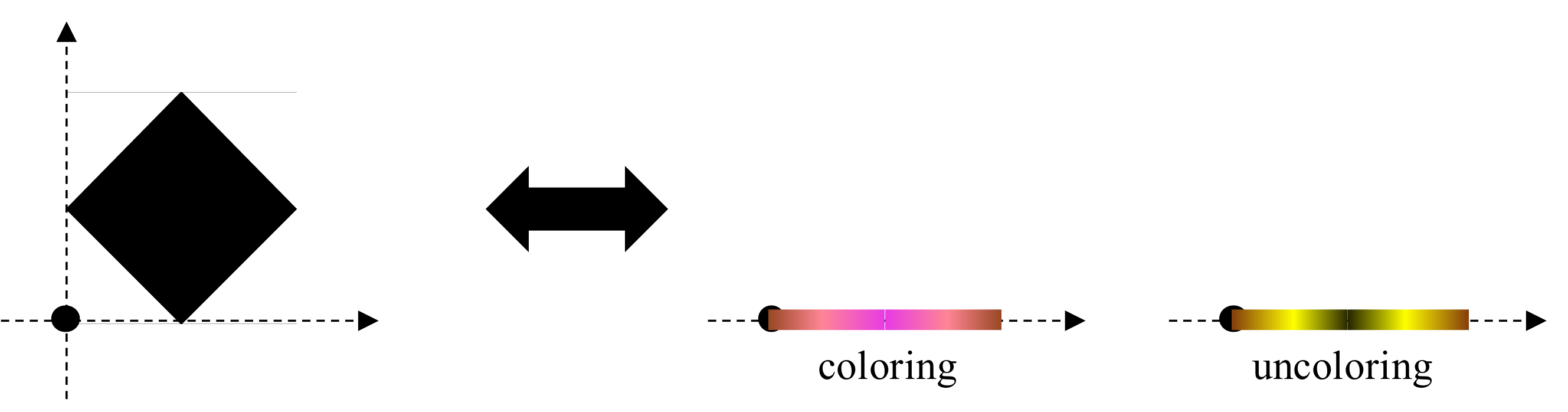}}}\hspace{5pt}
\caption{Color representation of truncated 2-cube.}\label{truncated_square_2d}
\end{figure}

Fig.~\ref{truncated_square_2d} shows the color representation of a truncated square. Note that the truncation for the top vertex is represented by the gradation of the coloring and for the side vertices is represented by uncolroing. This asymmetry is caused by the fact that the color representation of a square involves a projection of it onto one-dimensional Euclidean space.

\begin{figure}[H]
\centering
{%
\resizebox*{\textwidth}{!}{\includegraphics{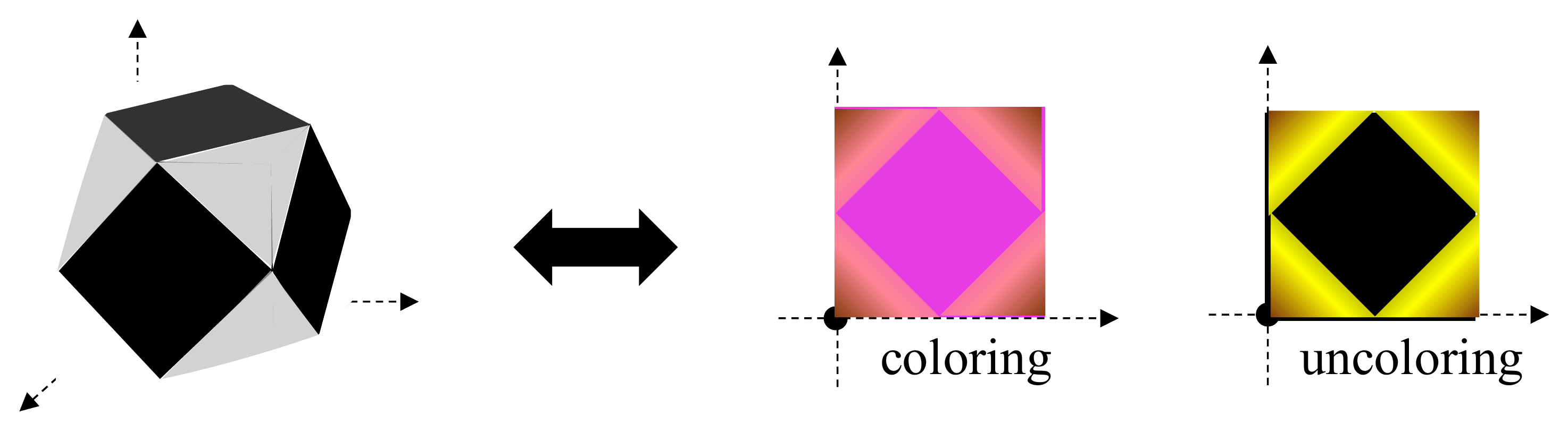}}}\hspace{5pt}
\caption{Color representation of truncated 3-cube.}\label{truncated_square_3d}
\end{figure}

\begin{figure}[H]
\centering
{%
\resizebox*{\textwidth}{!}{\includegraphics{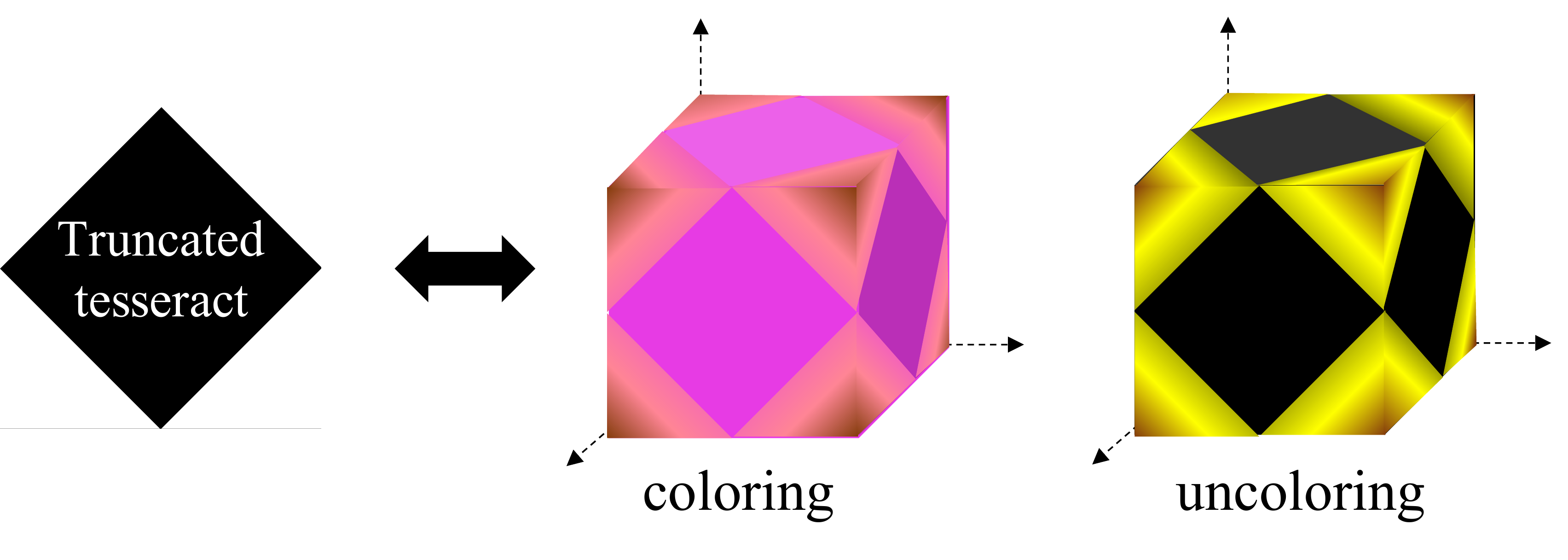}}}\hspace{5pt}
\caption{Color representation of truncated 4-cube.}\label{truncated_square_4d}
\end{figure}

Similarly, one can obtain the color representation of a truncated cube and tesseract by using the gradation of the coloring and the uncoloring (Figs.~\ref{truncated_square_3d} and \ref{truncated_square_4d}).

\begin{figure}[H]
\centering
{%
\resizebox*{\textwidth}{!}{\includegraphics{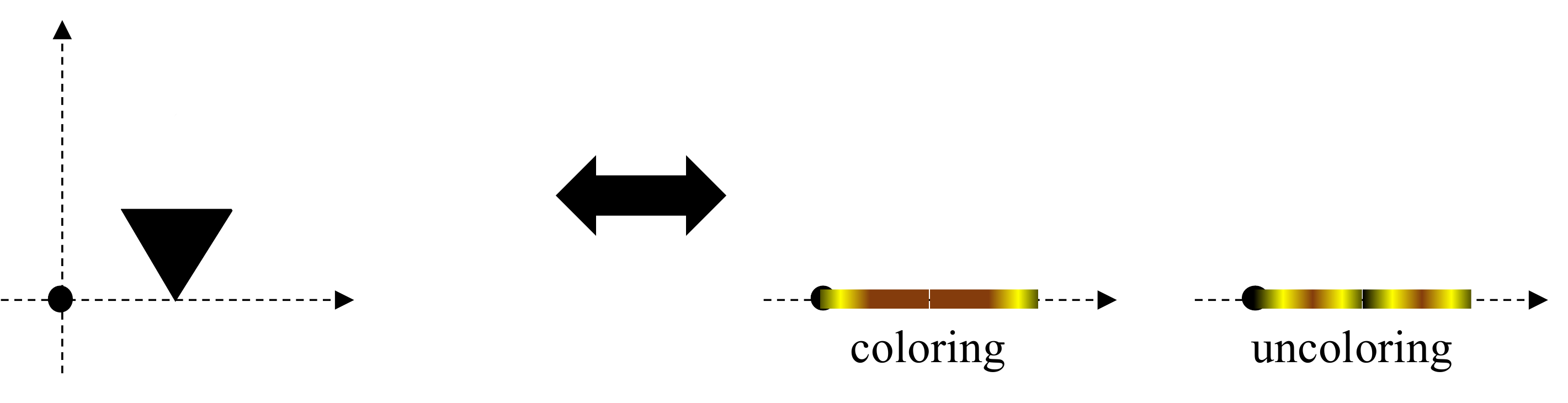}}}\hspace{5pt}
\caption{Color representation of truncated 2-simplex.}\label{truncated_simplex_2d}
\end{figure}

In addition, Fig.~\ref{truncated_simplex_2d} shows the representation of a truncated triangle.

\begin{figure}[H]
\centering
{%
\resizebox*{\textwidth}{!}{\includegraphics{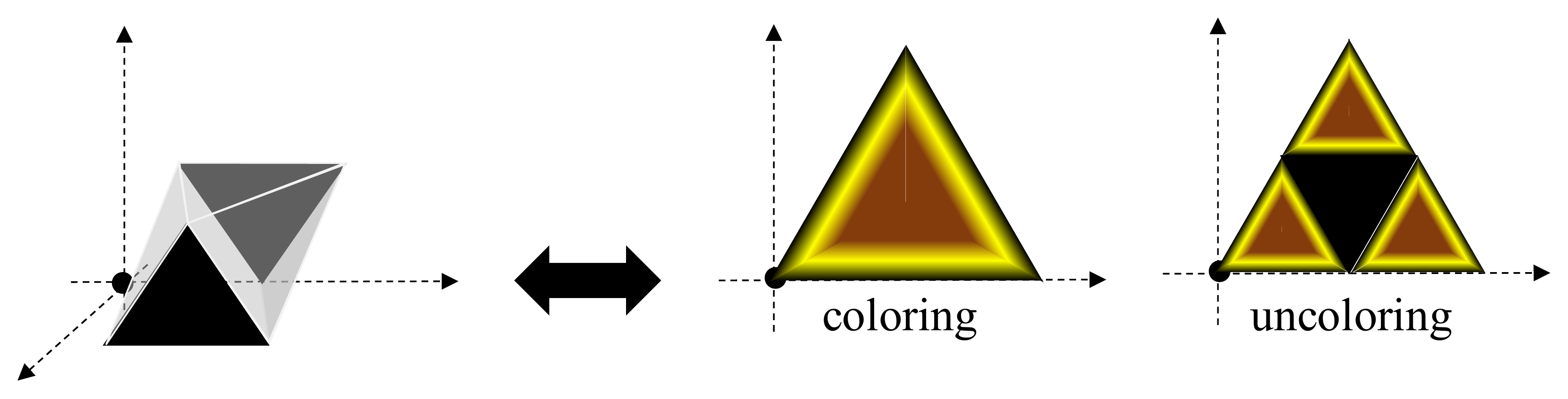}}}\hspace{5pt}
\caption{Color representation of truncated 3-simplex.}\label{truncated_simplex_3d}
\end{figure}

\begin{figure}[H]
\centering
{%
\resizebox*{\textwidth}{!}{\includegraphics{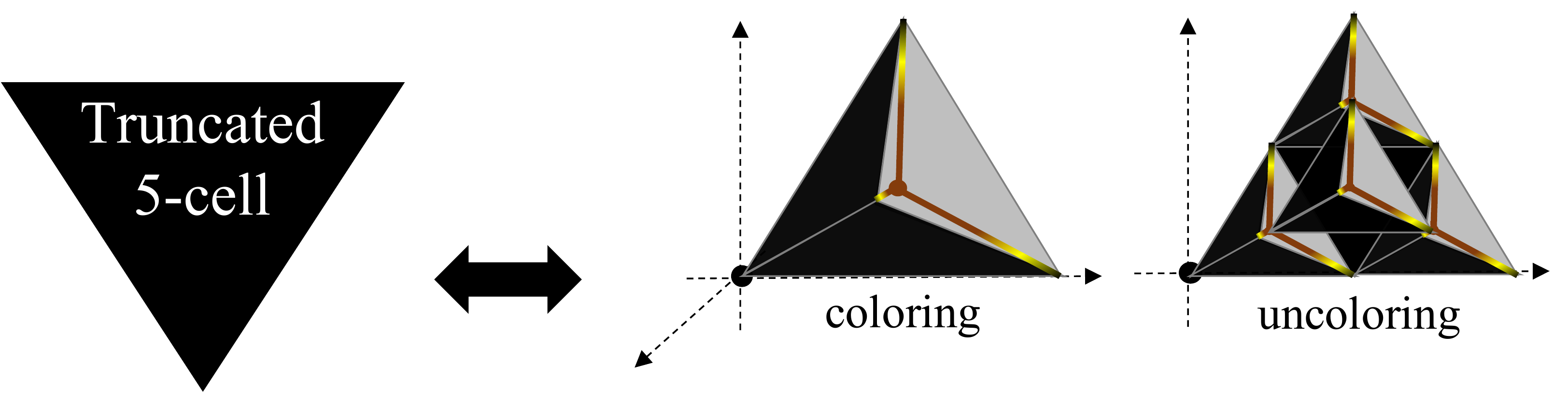}}}\hspace{5pt}
\caption{Color representation of truncated 4-simplex.}\label{truncated_simplex_4d}
\end{figure}

One can also obtain the color representation of a truncated tetrahedron (Fig.~\ref{truncated_simplex_3d}) and a truncated 5-cell (Fig.~\ref{truncated_simplex_4d}).

\newpage
\section{$1/n$-Coloring and Stellation}\label{1/n-Coloring and Stellation}
\begin{figure}[H]
\centering
{%
\resizebox*{8cm}{!}{\includegraphics{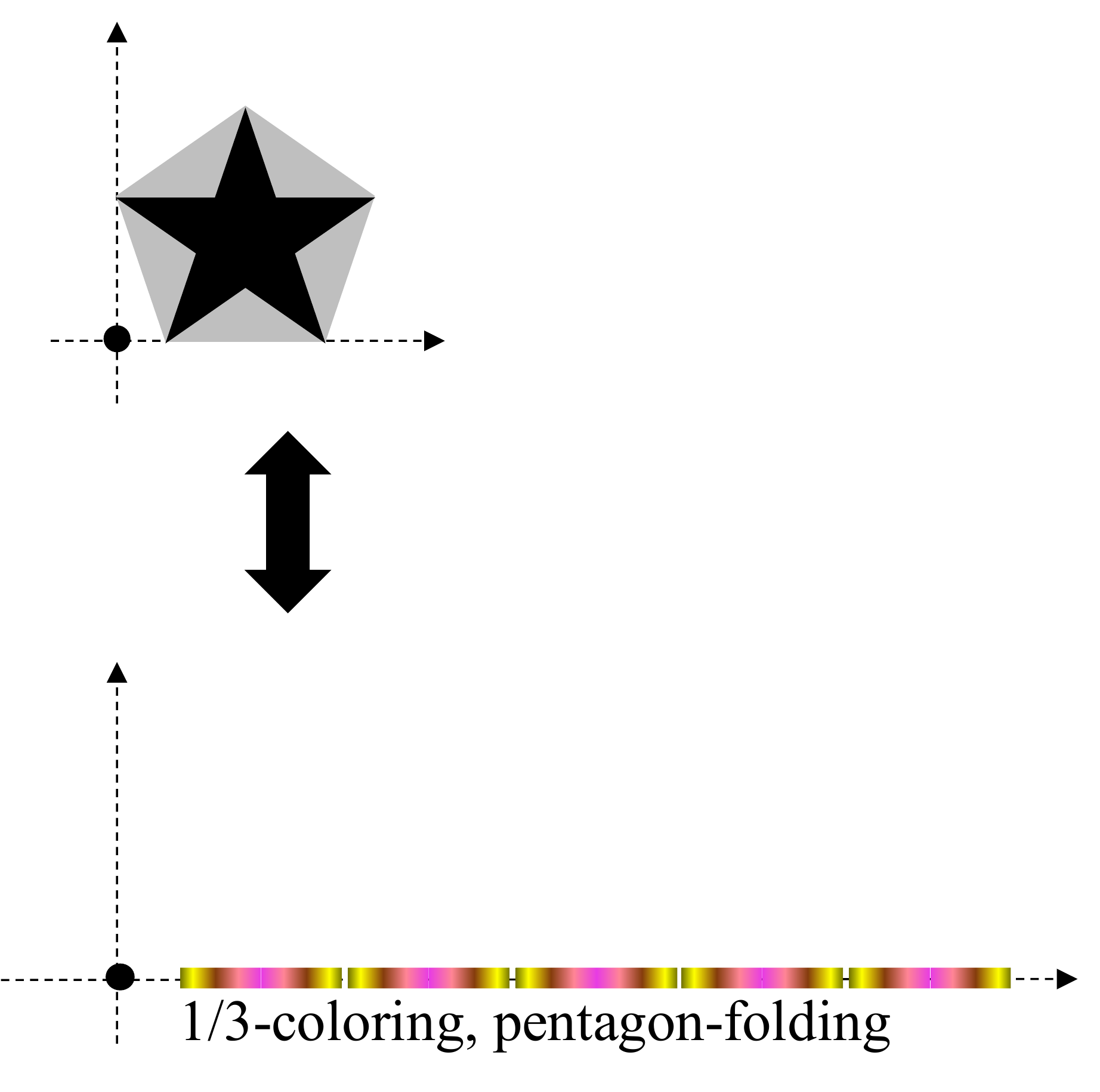}}}\hspace{5pt}
\caption{Color representation of a pentagram. The sum of the gray parts and the black part shows a pentagon that corresponds to the colored figure below if it is made by using usual coloring, not 1/3-coloring.}\label{star_5}
\end{figure}

\begin{figure}[H]
\centering
{%
\resizebox*{8cm}{!}{\includegraphics{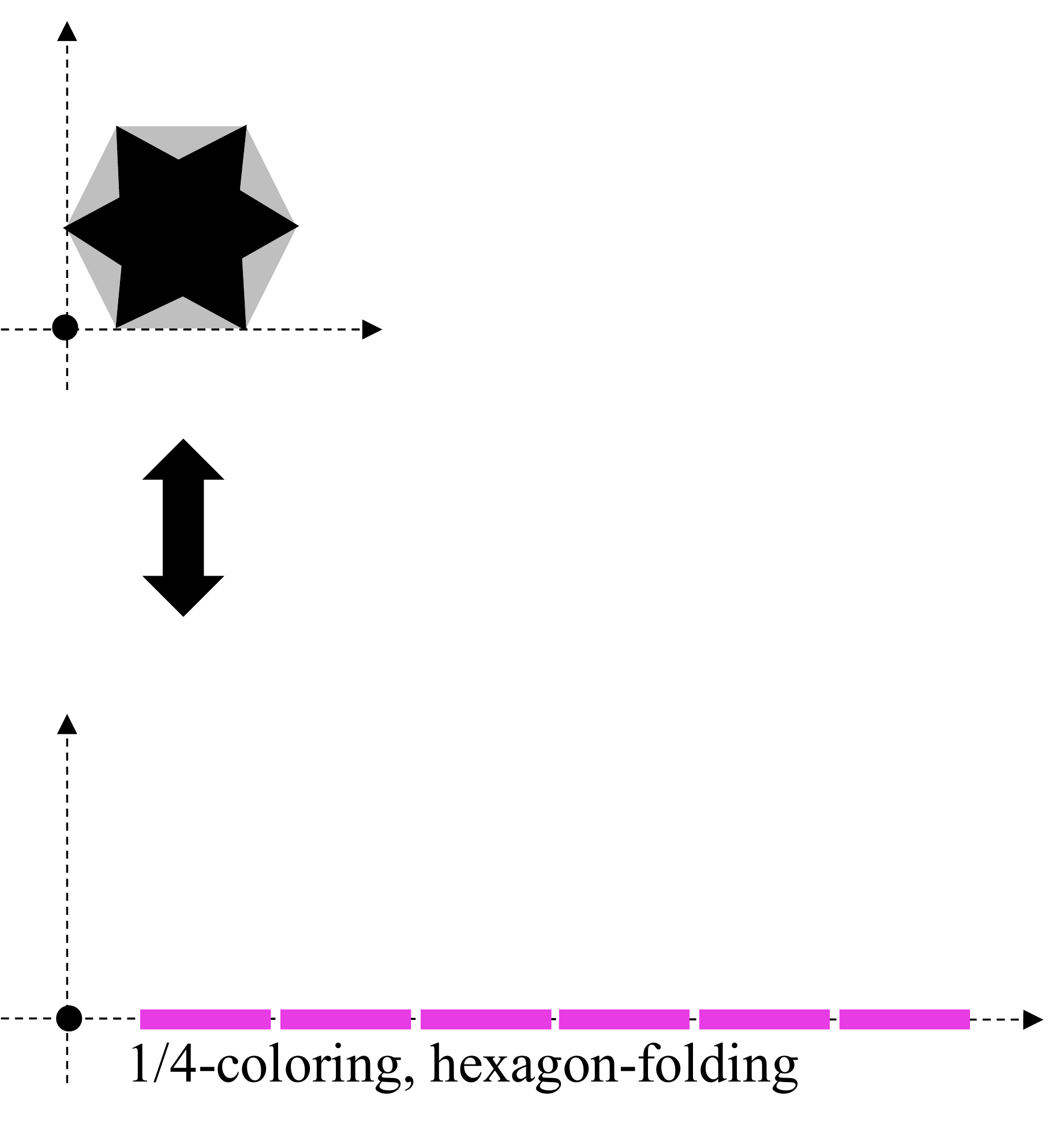}}}\hspace{5pt}
\caption{Color representation of a hexagram. The sum of the gray parts and the black part shows a hexagon that corresponds to the colored figure below if it is constructed by using usual coloring, not 1/4-coloring.}\label{star_6}
\end{figure}

\begin{figure}[H]
\centering
{%
\resizebox*{12cm}{!}{\includegraphics{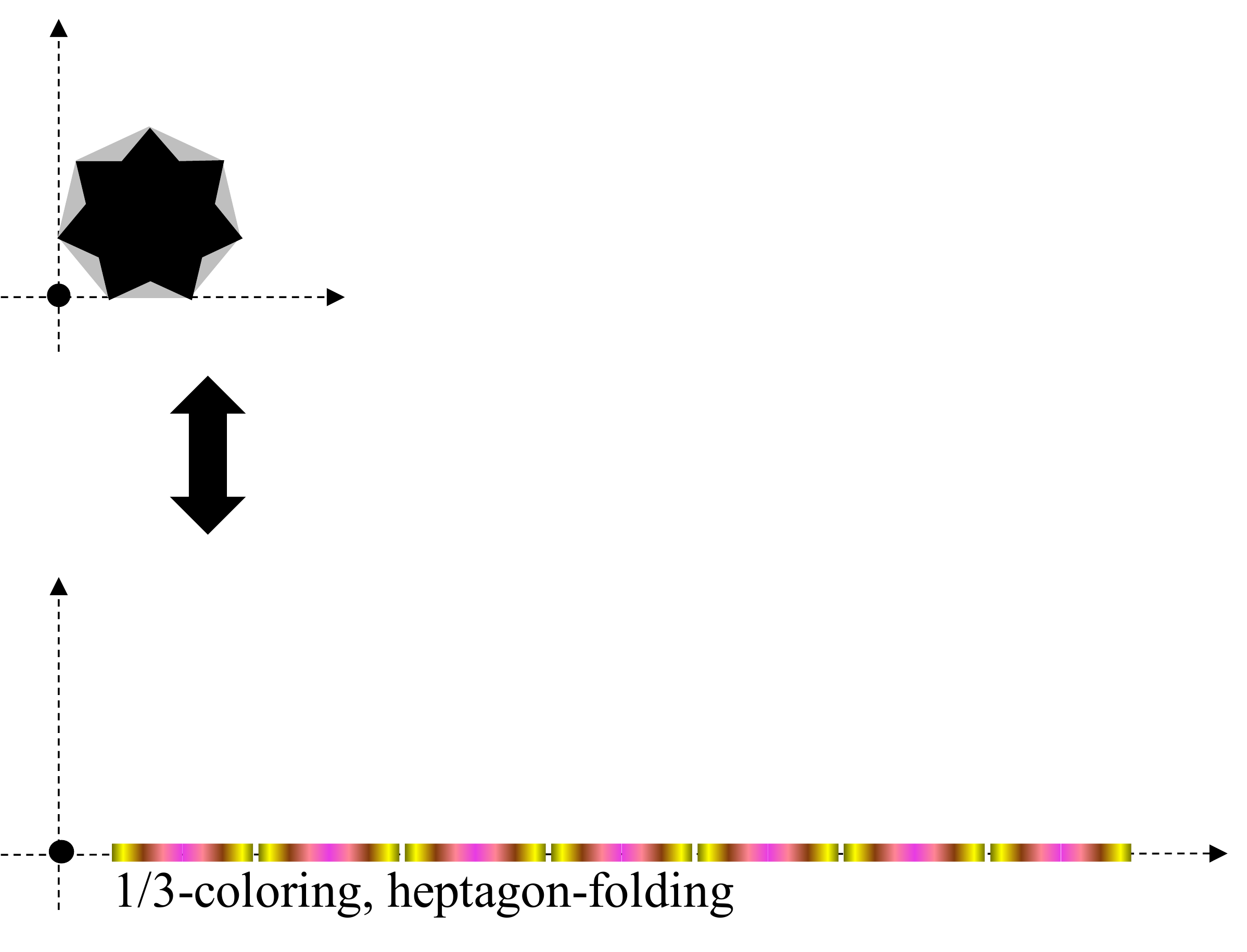}}}\hspace{5pt}
\caption{Color representation of a heptagram. The sum of the gray parts and the black part shows a heptagon that corresponds to the colored figure below if it is constructed by using usual coloring, not 1/3-coloring.}\label{star_7}
\end{figure}

\begin{figure}[H]
\centering
{%
\resizebox*{12cm}{!}{\includegraphics{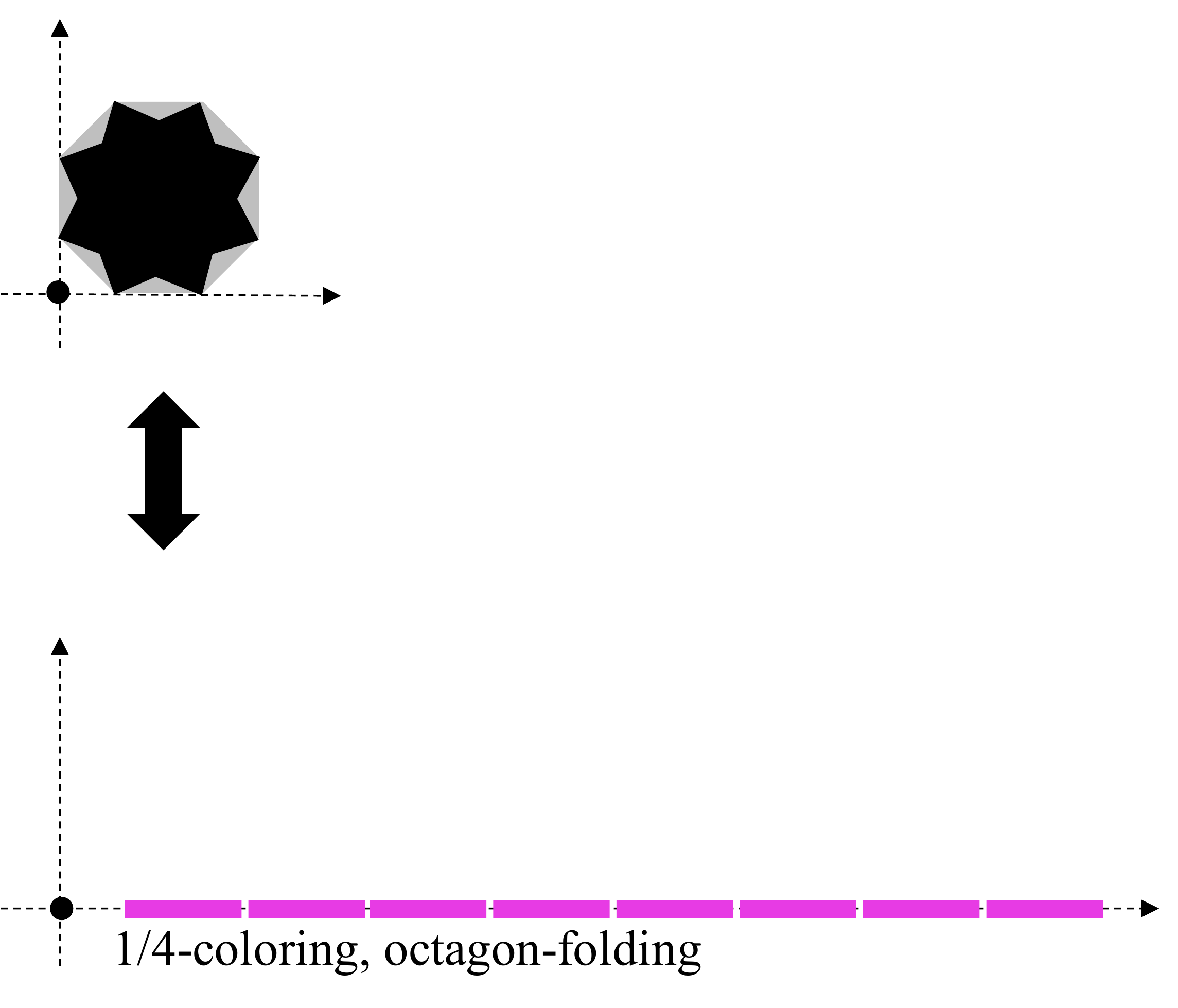}}}\hspace{5pt}
\caption{Color representation of an octagram. The sum of the gray parts and the black part shows a octagon that corresponds to the colored figure below if it is constructed by using usual coloring, not 1/4-coloring.}\label{star_8}
\end{figure}

In this section, we propose a way to represent a stellated polygon by using $1/n$-coloring introduced in Sec. \ref{Introduction}. Fig.~\ref{star_5} shows the unit pentagram by using five successive $1/3$-colored line segments with pentagon-folding. The maximum value for the color bar is set to $\sqrt{5+2\sqrt{5}}/2$. 
Similarly, Fig.~\ref{star_6} shows the unit hexagram by using six successive $1/4$-colored line segments with hexagon-folding. The maximum value for the color bar is set to $\sqrt{3}$. 
Fig.~\ref{star_7} shows the unit \{7/3\} heptagram by using seven successive $1/3$-colored line segments with heptagon-folding. The maximum value for the color bar is set to $\cot (\pi/14)/2$. 
Fig.~\ref{star_8} shows the unit \{8/3\} octagram by using eight successive $1/4$-colored line segments with octagon-folding. The maximum value for the color bar is set to $1+\sqrt{2}$.

\newpage
\section{Fractal}\label{Fractal}
In this section, we propose a way to represent some fractals by the coloring and uncoloring.

\begin{figure}[H]
\centering
{%
\resizebox*{8cm}{!}{\includegraphics{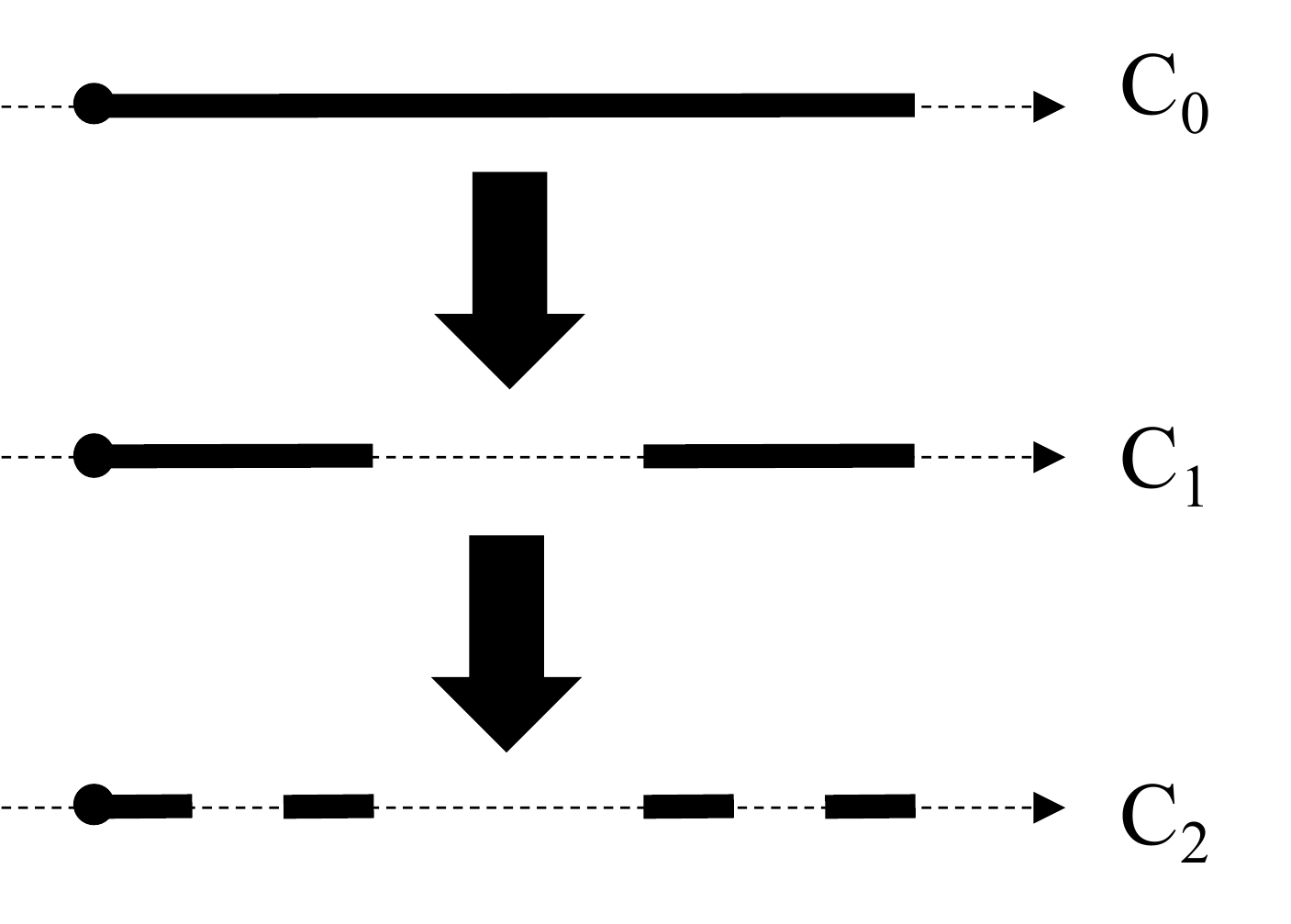}}}\hspace{5pt}
\caption{Iterations for creating Cantor set.}\label{Cantor}
\end{figure}

Cantor set is a fractal whose fractal dimension is $\log_32$. The Cantor set is considered as the limit $C_\infty$ of the following iterations: (1) $C_1$ is obtained by removing the center 1/3 from a line segment [0,1], (2) $C_2$ is obtained by removing the center 1/3 from the line segments obtained in $C_1$, (3) $C_3$ is obtained by removing the center 1/3 from the line segments obtained in $C_2$ (Fig.\ref{Cantor}).

\begin{figure}[H]
\centering
{%
\resizebox*{\textwidth}{!}{\includegraphics{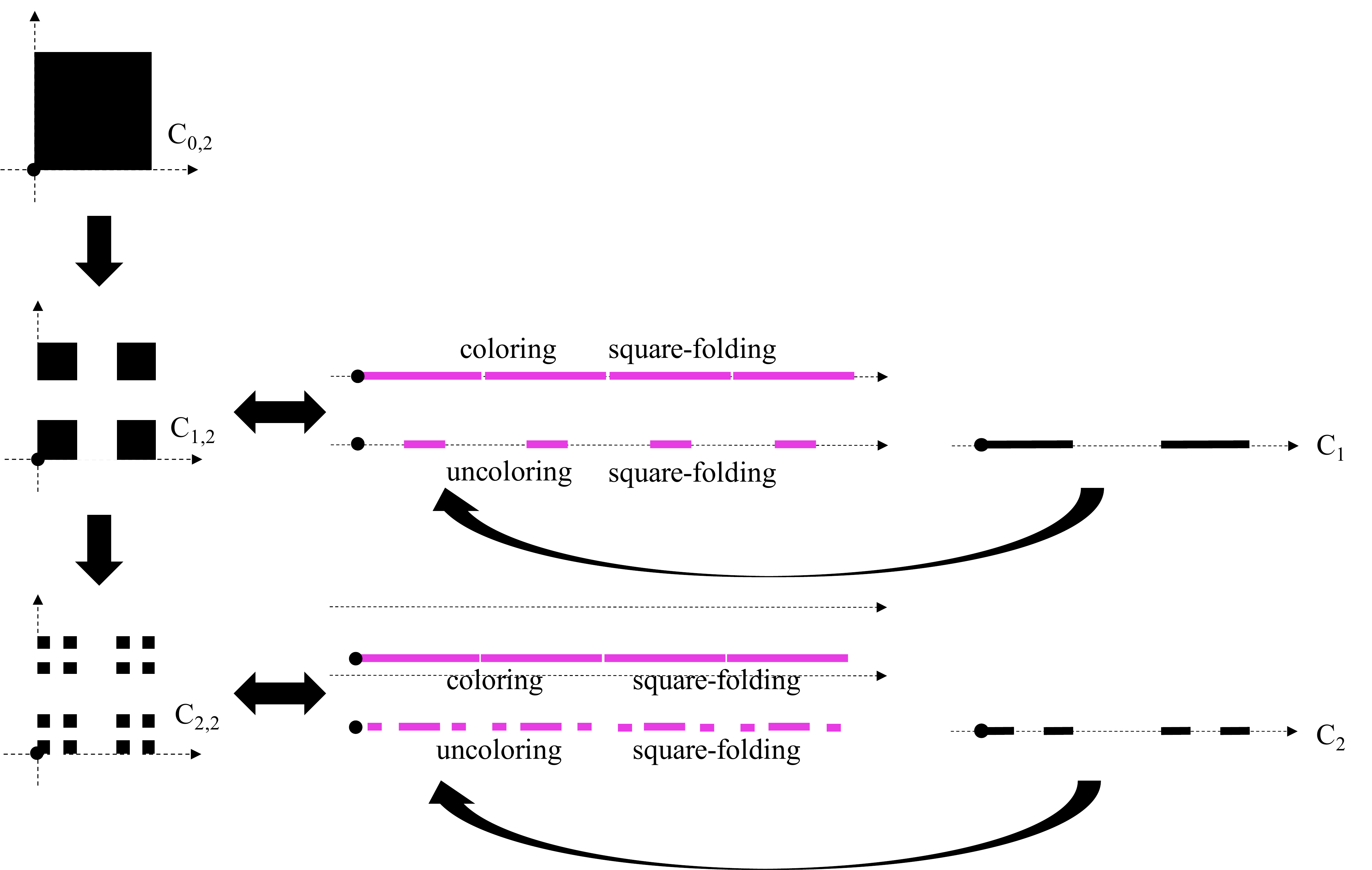}}}\hspace{5pt}
\caption{Iterations for creating 2-dimensional Cantor dust and their color representation.}\label{Cantor_dust_2d}
\end{figure}

\begin{figure}[H]
\centering
{%
\resizebox*{\textwidth}{!}{\includegraphics{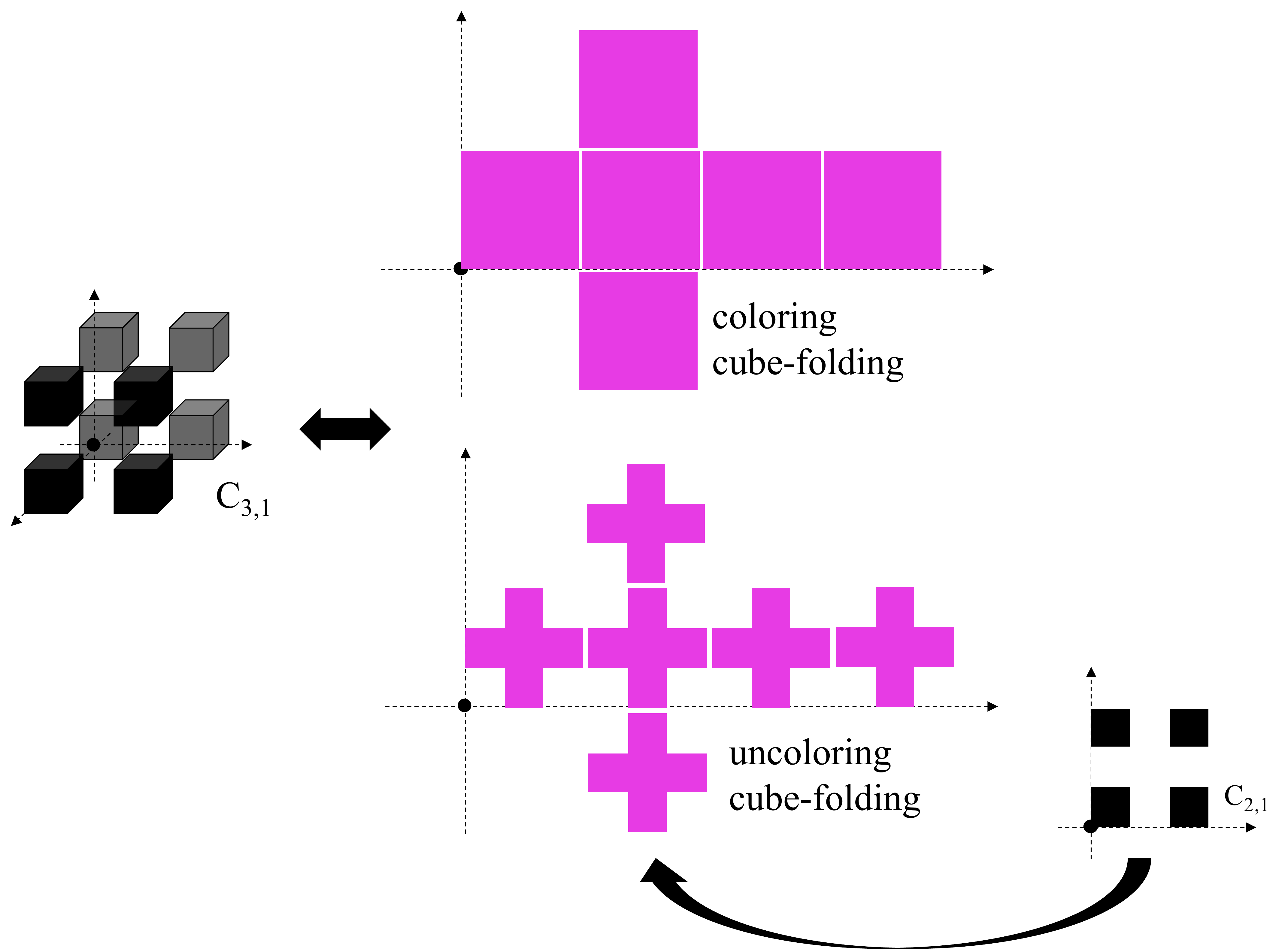}}}\hspace{5pt}
\caption{The first iteration for creating 3-dimensional Cantor dust and its color representation.}\label{Cantor_dust_3d}
\end{figure}

\begin{figure}[H]
\centering
{%
\resizebox*{\textwidth}{!}{\includegraphics{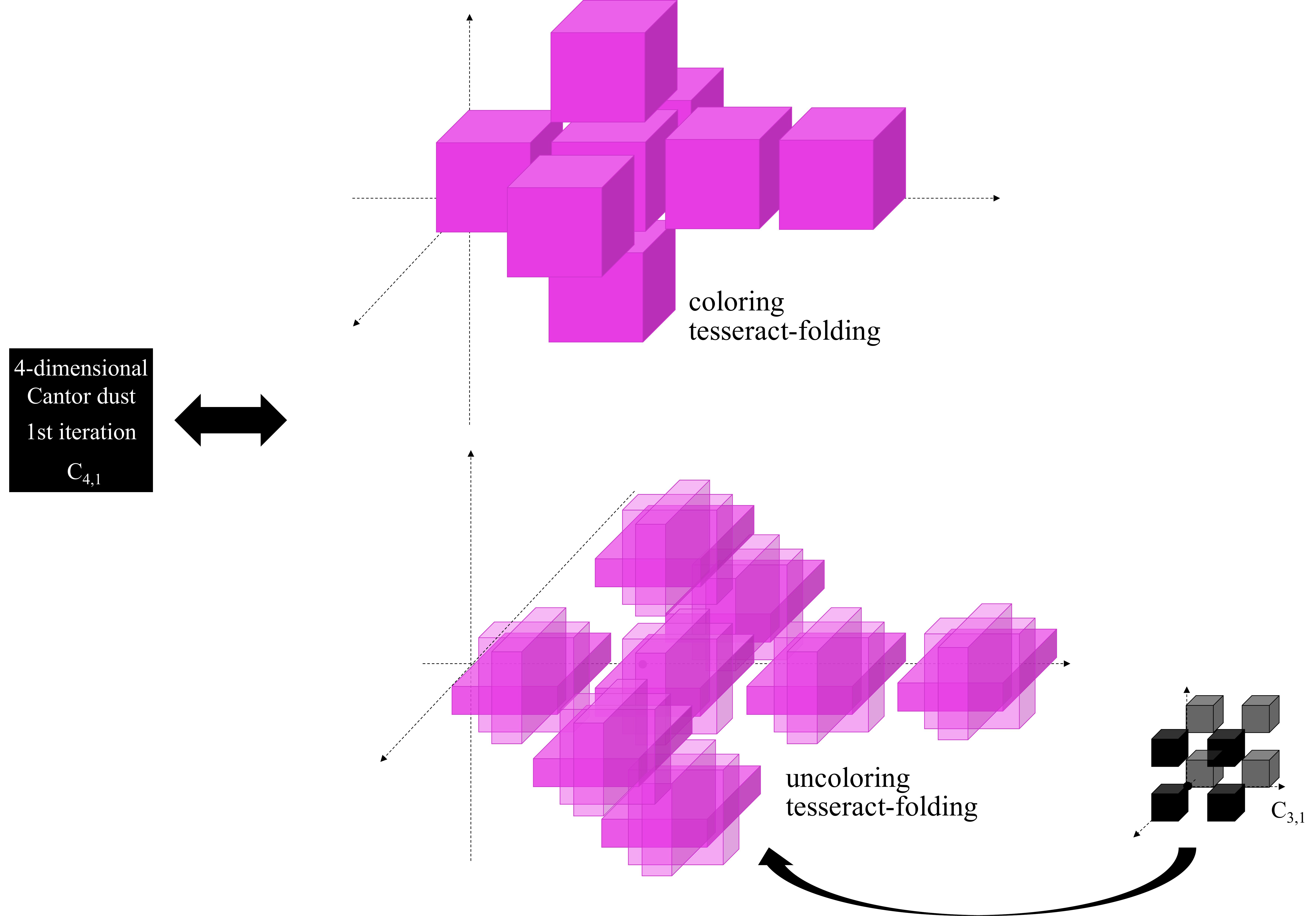}}}\hspace{5pt}
\caption{The first iteration for creating 4-dimensional Cantor dust and its color representation.}\label{Cantor_dust_4d}
\end{figure}

Then, we propose the color representation of the two-dimensional Cantor dust that is the Cartesian product of two Cantor sets by using the pink-uncoloring and square-folding. As with the Cantor set, the two-dimensional Cantor dust is considered as the limit $C_{\infty,2}$ of the iterations shown in Fig.~\ref{Cantor_dust_2d}. Then, we show how each iteration is represented by the iteration for Cantor set $C_n$ and various coloring techniques.  First, we replace $[0,1] \backslash C_n$ by its pink uncoloring version and making four copies and place them side by side (Fig.~\ref{Cantor_dust_2d}). Then, by combining a figure obtained by applying the square-folding to them and pink line segments, we have the figures shown in Fig.~\ref{Cantor_dust_2d}. Similarly, the three-dimensional Cantor dust can be represented by applying the pink uncoloring and the cube-folding to the six $[0,1]\times[0,1] \backslash C_{n,2}$ (Fig.~\ref{Cantor_dust_3d}). In addition, one obtains the color representation of the four-dimensional Cantor dust (Fig.~\ref{Cantor_dust_4d}).

\begin{figure}[H]
\centering
{%
\resizebox*{12cm}{!}{\includegraphics{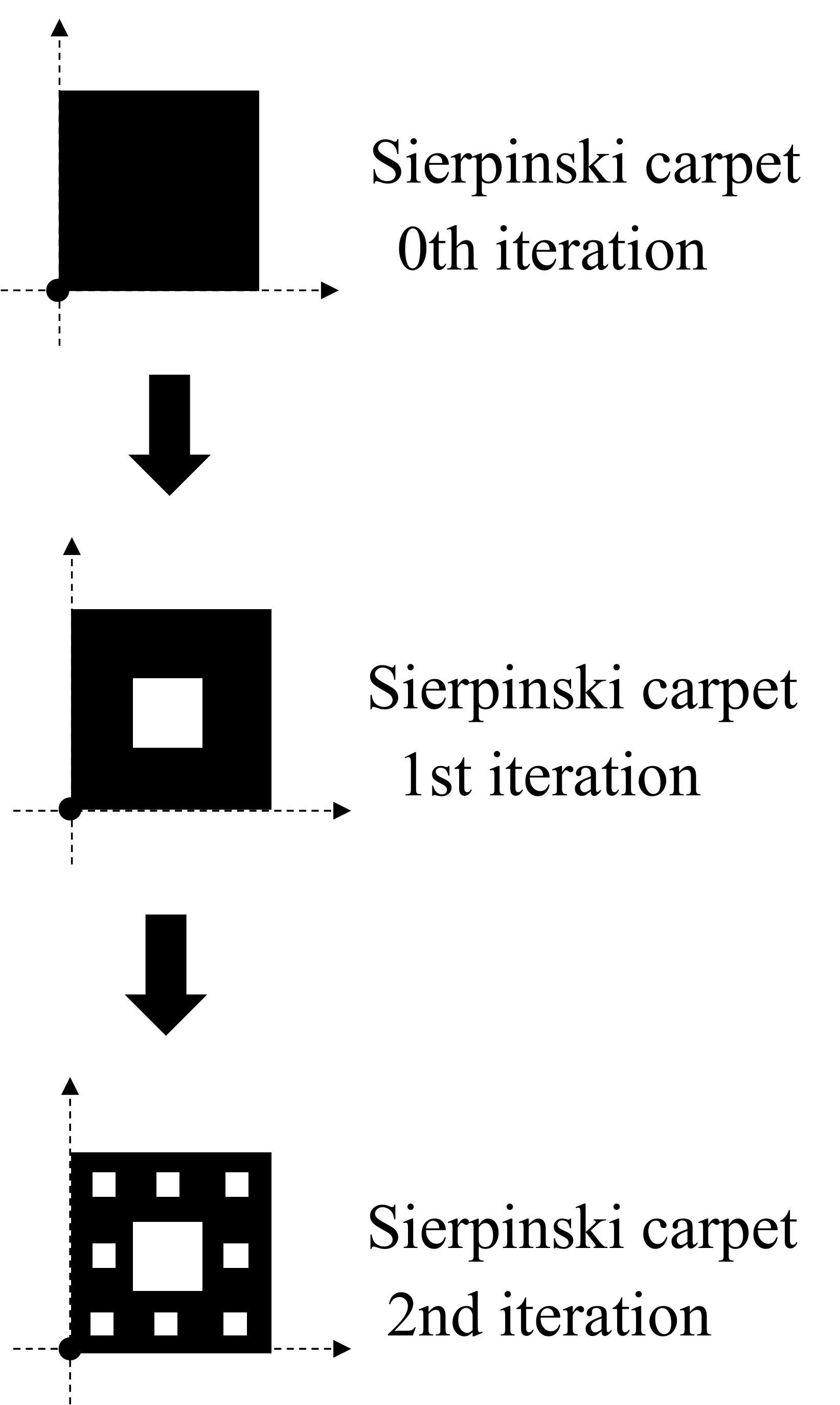}}}\hspace{5pt}
\caption{Iterations for creating Sierpinski carpet.}\label{Sierpinski}
\end{figure}

\begin{figure}[H]
\centering
{%
\resizebox*{\textwidth}{!}{\includegraphics{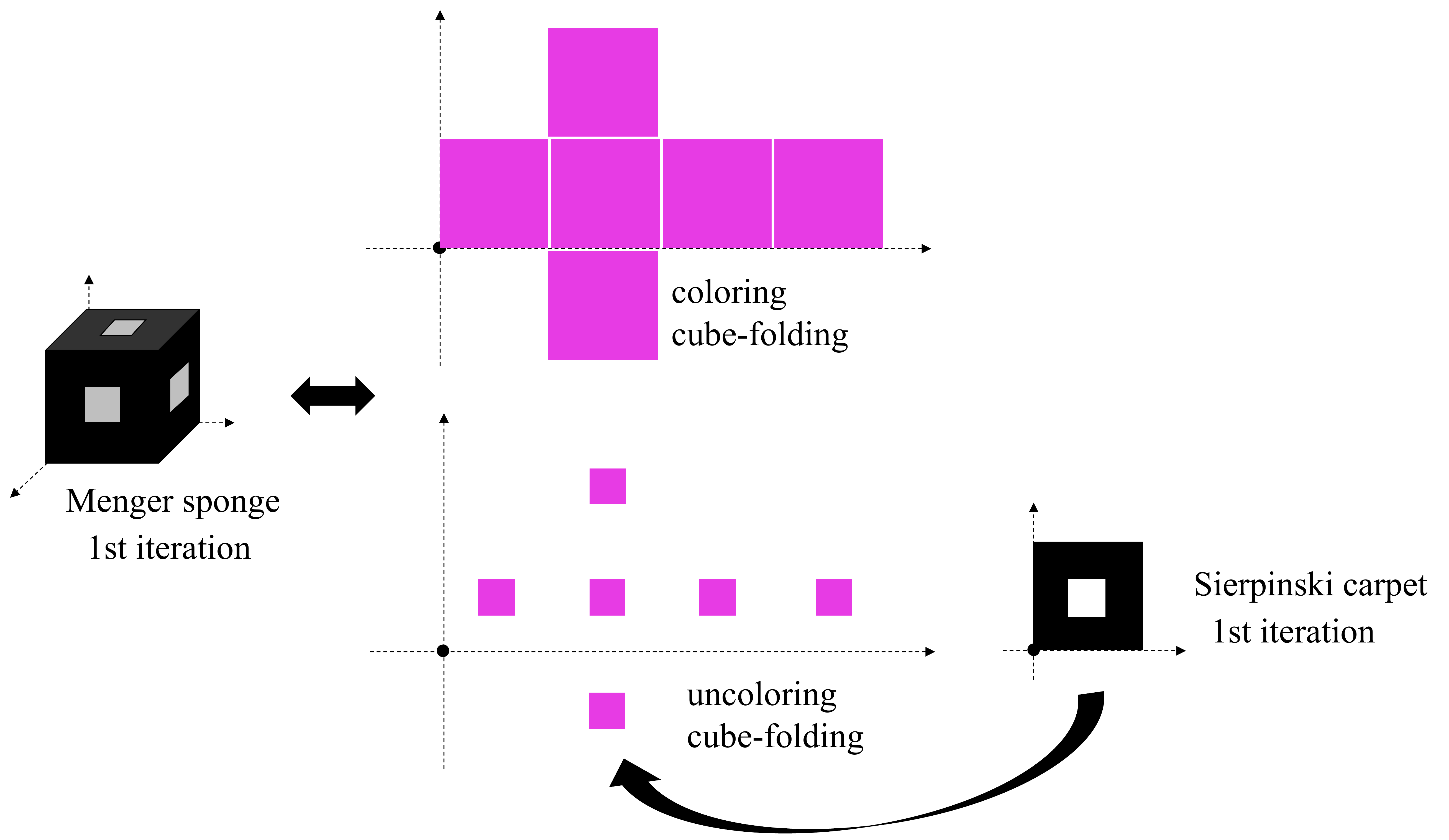}}}\hspace{5pt}
\caption{The first iteration for creating Menger sponge and its color representation.}\label{Menger}
\end{figure}

\begin{figure}[H]
\centering
{%
\resizebox*{12cm}{!}{\includegraphics{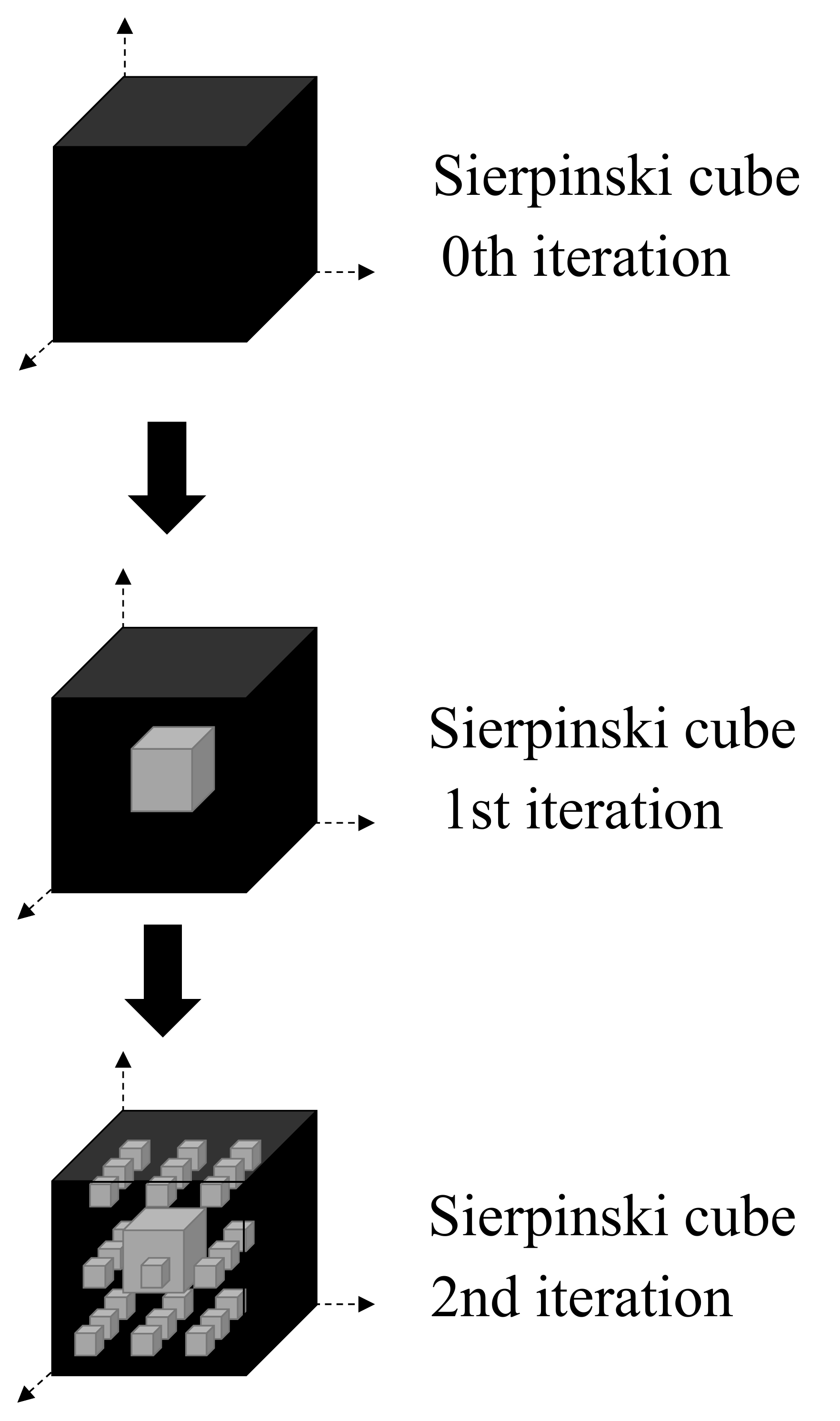}}}\hspace{5pt}
\caption{Iterations for creating Sierpinski cube where the gray parts represent holes.}\label{Sierpinski_cube}
\end{figure}

\begin{figure}[H]
\centering
{%
\resizebox*{\textwidth}{!}{\includegraphics{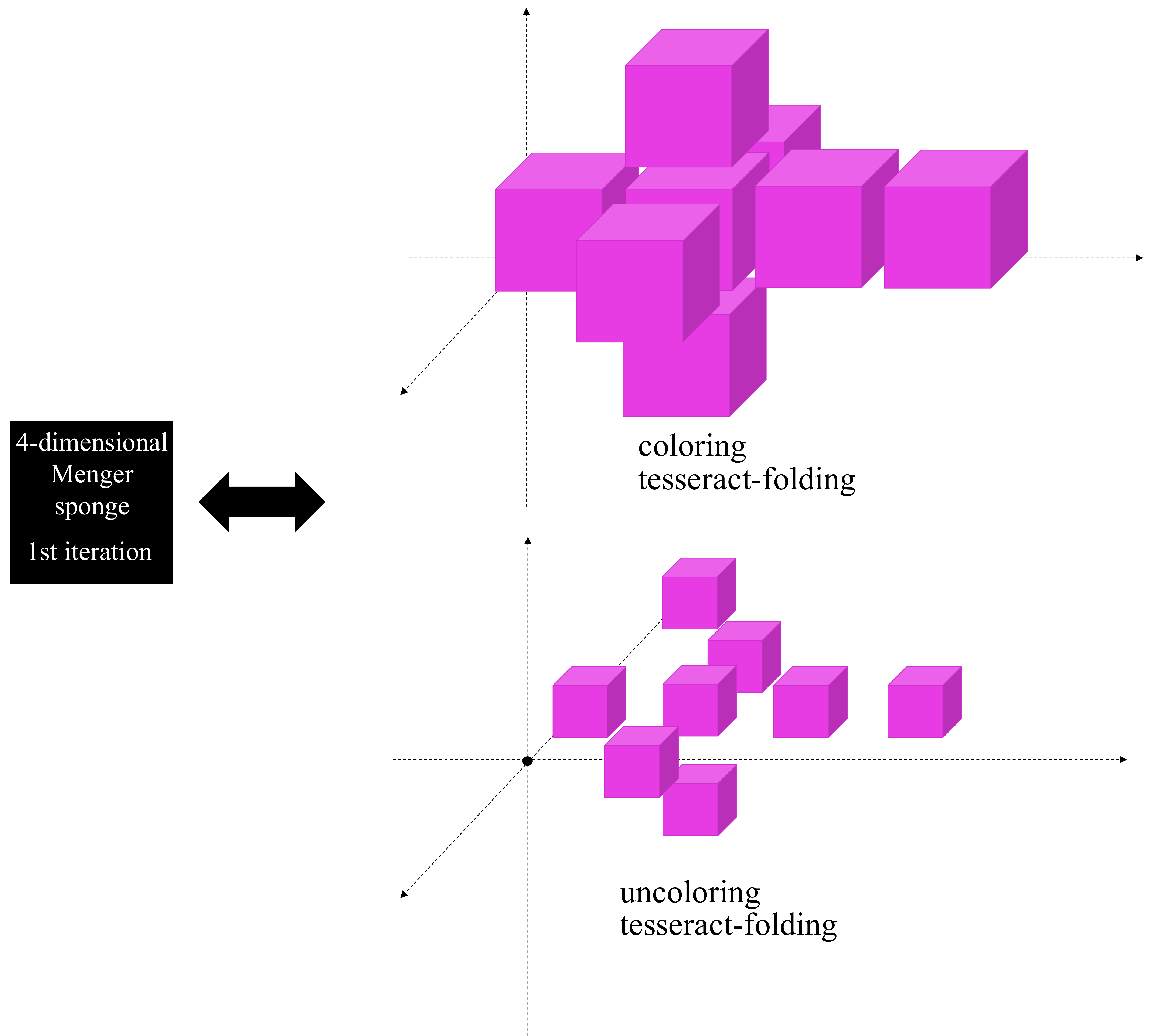}}}\hspace{5pt}
\caption{The first iteration for creating the 4-dimensional Menger sponge and its color representation.}\label{Menger_4d}
\end{figure}

Sierpinski carpet is a fractal considered as the limit of the iterations shown in Fig.~\ref{Sierpinski}. Menger sponge is the three-dimensional generalization and the color representation of the first iteration of the Menger sponge is shown in Fig.~\ref{Menger}. Also, one obtains color representation of the ``four-dimensional" version of the Menger sponge by using the cube instead of the square in each iteration of the usual Menger sponge and each iteration of the ``Siepinski cube" (Fig. \ref{Sierpinski_cube}) instead of each iteration of Sierpinski carpets in each iteration of the usual Menger sponge (see Fig.~\ref{Menger_4d} for the first iteration). The four-dimensional version of the Menger sponge whose fractal dimension is $\log_372 > 3$ and volume is infinite.

\newpage
\section{Conclusion} 
We proposed ways to represent geometrical objects by color. We introduced various actions to geometrical objects: coloring, uncoloring, $1/n$-coloring. First, we proposed ways to represent some regular polytopes (cube and simplex in four or fewer dimensions), the geometrical net of some regular polytopes (cube and simplex in five or fewer dimensions), and hyperprisms by using the coloring. Next, we introduced the concept of uncoloring and propose a way to represent truncated polytopes. Also, by introducing the concept of $1/n$-coloring, we proposed a way to represent stellated polytopes. Last, we represented the two-dimensional Cantor dust (three-dimensional Cantor dust) by the Cantor set (two-dimensional Cantor dust) and the Menger sponge by the Sierpinski carpet through the coloring and uncoloring, and propose the ``four-dimensional" Menger sponge whose volume is infinite.

In this paper, we have not dealt with some polytopes: orthoplex in $n$ dimensions and some regular polytopes in two, three, and four dimensions. Generally, it is difficult to deal with complex geometrical objects by color in the usual way because the resolution of the human eye is limited. In addition, the representation of geometrical objects by using color that we proposed loses their symmetry of them. For example, it is difficult to know the symmetry of a 4-dimensional cube by using its color representation of it. Once these issues are resolved, the color representation of geometrical objects will become even more important.


\end{document}